\documentclass[preprint,12pt]{elsarticle}
\usepackage{graphicx}
\usepackage{amsmath}
\usepackage{amssymb}
\usepackage{rotating}
\begin{document}

\begin{frontmatter}

%\title{How to unveil meaningful correlations through bifurcations: analysis of
%a mean field model of cortical activity}

\title{Metabifurcation analysis of a mean field model of the cortex}

\author{Federico Frascoli \corref{cor1}}\address{Brain Dynamics Research Unit\\Brain Sciences Institute
\\Swinburne University of Technology\\Hawthorn, Victoria\\AUSTRALIA}\ead{ffrascoli@swin.edu.au}
\author{Lennaert van Veen}\address{Faculty of Science\\University of Ontario Institute of Technology\\Oshawa, Ontario\\CANADA}
\author{Ingo Bojak}\address{Donders Institute for Brain, Cognition and Behaviour\\Centre for Neuroscience
\\Radboud University Nijmegen (Medical Centre)\\Nijmegen\\THE NETHERLANDS}
\author{David T J Liley}\address{Brain Dynamics Research Unit\\Brain Sciences Institute
\\Swinburne University of Technology\\Hawthorn, Victoria\\AUSTRALIA}
\cortext[cor1]{Corresponding author}

\begin{abstract}
% Background
Mean field models (MFMs) of cortical tissue incorporate
salient, average features of neural masses in order to
model activity at the population level, thereby linking
microscopic physiology to macroscopic observations, e.g.,
with the electroencephalogram (EEG).
One of the common aspects of MFM descriptions is the
presence of a high dimensional parameter space capturing
neurobiological attributes deemed relevant
to the brain dynamics of interest.

% Methodology/Principal Findings
We study the physiological parameter space of a MFM
of electrocortical activity and discover robust 
correlations between physiological attributes of the
model cortex and its dynamical features. These correlations are
revealed by the study of bifurcation plots, which show that the model 
responses to changes in inhibition belong to two archetypal categories 
or ``families''. After investigating and characterizing them in depth,
we discuss their essential differences in terms of four important
aspects: power responses with respect to the modeled action of anesthetics, reaction
to exogenous stimuli such as thalamic input, 
distribution of model parameters and
oscillatory repertoires when inhibition is enhanced.
Furthermore, while the complexity of sustained periodic orbits differs 
significantly between families, we are able to show how metamorphoses
between the families can be  brought about by exogenous stimuli.

% Conclusions/Significance
We here unveil links between measurable physiological
attributes of the brain and 
dynamical patterns that are not accessible by linear methods. 
They instead emerge when the nonlinear structure of parameter space is
partitioned according to bifurcation responses.
This partitioning cannot be achieved by the investigation of only a small
number of parameter sets and is instead the result of an automated bifurcation analysis
of a representative sample of 73,454 physiologically admissible parameter sets.
Our approach generalizes straightforwardly and is well suited to
probing the dynamics of other models with
large and complex parameter spaces.
\end{abstract}


\end{frontmatter}

\section{Introduction}

Understanding brain function is one of the longstanding and unresolved quests of science.
Despite the efforts and contributions of a wide range of scientists from
different disciplines, we are still very far from a comprehensive account of 
what the brain does, and how and why it does it. The formidable hardness of this problem 
is embodied in the complexity of the structure and organization of this organ,
which supports interactions over a wide range of temporal and spatial scales.
One of the scientific approaches that aims to shed some light on the complexity of brain activity is
the use of realistic mathematical models and computer simulations: 
constrained by known physiology and anatomy, these theories try to
recreate and understand the patterns emerging in the electrical activity of the cortex.
In particular, so-called mean field models (MFMs) \cite{ref:JH96,ref:RRW97,ref:LCD02} 
try to capture the results of neuronal interactions through spatial and/or temporal averaging,
which then describe the population activity of regions of cortex in a parsimonious manner.
Accompanied by a sufficiently flexible parametrization, MFMs
take into account existing uncertainties in cortical anatomy and physiology
in an effective way.

Despite the apparent biological simplicity of MFMs, 
such theories  produce physiologically and behaviorally
plausible dynamics. Since their birth forty years ago, 
they have been employed successfully to investigate
a large number of cognitively and clinically relevant phenomena, such as the action
pharmacological agents have on the activity of cortical tissue 
\cite{ref:WC72,ref:WC73,ref:LHSZ73,ref:Fre75,ref:LvRB+76,ref:ZKM78},
the origin of evoked potentials \cite{ref:JR95,ref:RRW02},
the sleep-wake cycle and circadian rhythm \cite{ref:SSS+05,ref:PR07}, 
the emergence of gamma band
oscillations and its relationship to cognition \cite{ref:Wri97,ref:RWR00,ref:BL07}
and the onset and characteristics of epileptic
seizures \cite{ref:WBBC00,ref:RRR02,ref:LBK+03,ref:KKS05,ref:LB05}; 
c.f.\ Ref.~\cite{ref:DJR+08} for a recent review. 
We use here Liley's MFM \cite{ref:LCD02}, which is
capable of reproducing the main spectral features of spontaneous
(i.e., not stimulus-locked) EEG, in particular the ubiquitous alpha
rhythm \cite{ref:LAWA99}. Firstly, this model is biologically constituted: 
all state variables and parameters can be constrained on the basis
of existing anatomical and physiological measurements in the literature. 
Secondly, it supports a rich repertoire of behaviors both 
physiologically relevant and dynamically interesting. For example, 
parametrically widespread, robust chaotic activity of
various origins has been found \cite{ref:DLC01,ref:vVL06,ref:FvVB+08}, 
and multistability, i.e., the presence of various coexisting
dynamical regimes, has been demonstrated. Multistability has been
speculated to correspond neurobiologically to the formation of
memories \cite{ref:DFCL09}. 

In the spirit of the dynamical approach \cite{ref:Fre01}, in this paper we link
nonlinear electrical activity and neurobiologically significant attributes
of cortex. 
To this end, we consider a representative sample of parameter sets that
have previously been found to generate physiologically plausible behavior \cite{ref:BL05}.
This sample contains 73,454 sets and for each of these we computed the bifurcation
plots when varying two parameters related to inhibition,
i.e., the qualitative changes of recurrent  activity patterns, when neuronal inhibition
is altered.  It turns out that we can sort the sets into two
distinct ``families'' of dynamical behavior. These families are found 
to correlate with EEG signal power and responses to anesthetics,
whereas family membership is determined by specific neurobiological parameters.

The paper is divided into three main sections. First, we briefly introduce 
Liley's ordinary differential Equations (ODEs). 
We next discuss the theoretical and numerical tools employed in some detail, in particular how
one arrives at a systematic
bifurcation analysis procedure to show the qualitative changes of the model solutions
when inhibition is varied. 
Then, we present the main distinctive features of the families, such as 
their responses to the simulated induction of anesthetics and their correlations
with model parameters of interest. We are also able to show how exogenous agents, i.e., 
input from the thalamus to the cortex, can induce dramatic changes in those patterns
and stimulate transitions from one type of family to the other.
That, in particular, provides a compelling example for the modulation thalamus is thought to
exert on the cortex \cite{ref:SG06}. All these relations cannot be discovered with standard
linear or nonlinear analyses of the physiological parameter space, and represent the main
result of this work. A discussion of open problems concludes the paper.

\section{Neuronal mean field equations}

Liley's MFM aims to provide a mathematically and physiologically parsimonious description of average neuronal activity in a
human cortex, with spatially coarse-grained but temporally precise dynamics.
One excitatory and one inhibitory neuronal population, respectively,
is considered per macrocolumn, which is a barrel-shaped region 
of approximately $0.5-3~\mathrm{mm}$ diameter comprising the whole thickness
of cortex (thus $\approx 3-4~\mathrm{mm}$ deep). 
Cortical activity is locally described by the mean soma membrane potentials of
the excitatory ($h_e$) and inhibitory ($h_i$) neuronal populations, along with four mean
synaptic inputs $I_{ee}$, $I_{ie}$, $I_{ei}$, and $I_{ii}$. These inputs convey
the reciprocal interaction between neuronal populations, where double subscripts
indicate first source then target (each either excitatory $e$ or inhibitory $i$). 
The connection with measurements is through $h_e$, which is linearly related to
the EEG signal \cite{ref:Nun81}. 
Lumped neuron populations are modeled as passive $RC$ compartments, into which all
synaptically induced ionic currents terminate. According to population types $(j,k) = e, i$,
synaptic activity drives the mean soma membrane potentials from their resting values. 
The equations for $h_e$ and $h_i$ are given by
\begin{eqnarray}
\tau_{e}\frac{d h_{e}}{dt} &=& h_e^r - h_e(t) + 
\frac{h^{eq}_{ee}-h_{e}(t)}{|h^{eq}_{ee}-h^r_e|}I_{ee}(t)
+\frac{h^{eq}_{ie}-h_{e}(t)}{|h^{eq}_{ie}-h^r_e|}I_{ie}(t)\;,\\
\label{eq:liley1}
\tau_{i}\frac{dh_{i}}{dt} &=& h_i^r - h_i(t) + 
\frac{h^{eq}_{ei}-h_{i}(t)}{|h^{eq}_{ei}-h^r_i|}I_{ei}(t)
+\frac{h^{eq}_{ii}-h_{i}(t)}{|h^{eq}_{ii}-h^r_i|}I_{ii}(t)\;,
\label{eq:liley2}
\end{eqnarray}
where $h^r_e$ and $h_i^r$ are mean resting potentials, and $\tau_{e}$ and
$\tau_{i}$ are the membrane time constants of the respective neuronal populations.  The
reversal potentials of the transmembrane ionic fluxes mediating excitation and inhibition
are given by $h^{eq}_{ek}$ and $h^{eq}_{ik}$, respectively. Note that
the synaptic inputs are weighted with $+1$ (excitatory $I_{ek}$) and $-1$ 
(inhibitory $I_{ik}$) at the resting potential of the respective excitatory or inhibitory
neuronal population, and that these weights then vary
linearly with voltage.

\begin{figure}[htbf]
\begin{center}
\includegraphics[width=\textwidth]{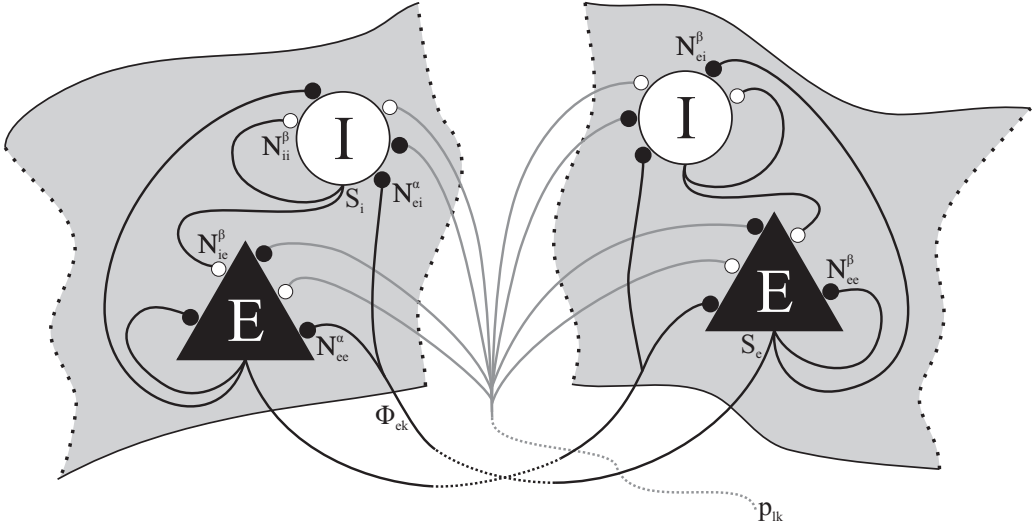}
\end{center}
\caption{{\bf Architecture of Liley's mean field model.} 
Two separate model macrocolumns are shown,
each containing one excitatory and one inhibitory neuronal population.
Note that long-range connections are exclusively excitatory and that
self-couplings correspond to connections of neurons of the same type
{\em within} the local populations.}
\label{fig:scheme}
\end{figure}

The mean synaptic inputs describe the postsynaptic 
activation of ionotropic neurotransmitter receptors by presynaptic action potentials, 
arising from the collective activity of neurons both nearby and distant. 
The time course of such activity, based on well-established experimental
data \cite{ref:Tuc88}, is modeled by a critically damped oscillator 
driven by the mean rate of incoming excitatory or inhibitory axonal pulses.
We thus have, for $k = e, i$:
\begin{eqnarray}
\left(\frac{d}{dt} + \gamma_{ek}\right)^2 I_{ek}(t) & = & 
\Gamma_{ek}\gamma_{ek}e\left\{N^\beta_{ek}S_{e}\left[h_{e}(t)\right] + p_{ek}(t) + \phi_{ek}(t)\right\}\;,
\label{eq:liley3} \\
\left(\frac{d}{dt} + \gamma_{ik}\right)^2 I_{ik}(t) & = & 
\Gamma_{ik}\gamma_{ik}e\left\{N^\beta_{ik}S_{i}\left[h_{i}(t)\right] + p_{ik}(t)\right\} \;,
\label{eq:liley4}
\end{eqnarray}
where the terms in curly brackets correspond to sources of the axonal pulses
from three origins:
local, i.e., in the same macrocolumn of the cortex $N^\beta_{lk}S_{l}$, 
arriving through long-range, excitatory cortico-cortical connections from other macrocolumns
$\phi_{ek}$, and extracortical, i.e., primarily of thalamic origin $p_{lk}$.
For subsequent simplicity we assume the absence of any extracortical inhibitory input, i.e. $p_{ik}\equiv 0$.
$N^\beta_{lk}$ quantifies the strength of anatomical population connectivity. 
The maximal postsynaptic potential (PSP) amplitude $\Gamma_{lk}$ occurs in the
target population $k=e,i$ at time $1/\gamma_{lk}$ after the arrival of the 
presynaptic spike from the source population $l=e,i$. 
A schematic illustration of the architecture of interactions in the
Liley model can be found in Fig.~\ref{fig:scheme}.

Local mean soma potentials $h_k$ are nonlinearly transformed
to mean neuronal population firing rates with a sigmoidal function
\begin{equation}
S_{k}\left[h_{k}(t)\right] = S_k^{max}
\left\{1+\exp\left[-\sqrt{2}\frac{h_k(t)-\mu_k}{\sigma_k}\right]\right\}^{-1}\;, 
\end{equation}
where $\mu_k$ and $\sigma_k$ indicate the firing thresholds and
their associated standard deviations for the respective neural population.
The chosen form is a computationally convenient approximation to the error function,
which results when one assumes a Gaussian distribution of firing thresholds.

The final part of the model concerns the propagation of excitatory
axonal activity in long-range fibers. In the full formulation, the following
inhomogeneous, two-dimensional damped wave equation is used:
\begin{equation}
\left[\left(\frac{\partial}{\partial t}+v\Lambda\right)^2
-\frac{3}{2}v^2\nabla^2\right]\phi_{ek}\left(\vec{x},t\right)
=N^\alpha_{ek}v^2\Lambda^2 S_{e}\left[h_{e}\left(\vec{x},t\right)\right]\;. \label{eq:liley5}
\end{equation}
Two core assumptions form the basis
of this equation: first, action potentials traverse much larger distances in the cortico-cortical fibres than
the local (intra-columnar) axonal fibres, and thus their conduction delays can no 
longer be considered negligible. Second, cortical areas further away share fewer long range connections.
Given the average number of long range excitatory connections onto a local neuron $N^{\alpha}_{ek}$,
experiments suggest that the fraction thereof that stems from distant neurons decays approximately
exponentially with distance, with a scale constant $\Lambda$. Activity is then modeled as spreading 
omnidirectionally with constant conduction speed $v$. We have also made the additional assumption that
long-distance propagation dynamics is not differentiated according to excitatory and inhibitory targets, i.e.,
 $\Lambda_{ee}=\Lambda_{ei}=\Lambda$ and $v_{ee}=v_{ei}=v$. 
 
In the following, we force the term $\nabla^2\phi_{ek} \rightarrow 0$ and thus consider only
bulk dynamics for our model cortex. Furthermore, we neglect the spatially inhomogeneous effects 
of long-range fibre systems. Thus, Eq.~(\ref{eq:liley5}) now becomes
\begin{equation}
\left(\frac{d}{dt}+\omega\right)^2 \phi_{ek}(t)
=N^\alpha_{ek}\omega^2 S_{e}\left[h_{e}(t)\right]\;, \label{eq:liley6}
\end{equation}
with $\omega=v\Lambda$.
We see that Eq.~(\ref{eq:liley6}) has turned into an inhomogeneous critically damped
oscillator. This is in contrast to simulating an unconnected local population only,
which would involve dropping Eq.~(\ref{eq:liley5}) altogether. However, it is clear that the steady
state solution to Eq.~(\ref{eq:liley6}), $\phi_{ek}(t)=N^\alpha_{ek}S_{e}\left[h_{e}(t)\right]$,
simply changes $N^\beta_{ek}\rightarrow N^\alpha_{ek}+N^\beta_{ek}$ in Eq.~(\ref{eq:liley3}). Thus,
the steady state of bulk oscillations corresponds to local activity with increased connectivity.

The fourteen
ODEs of Eqs.~(\ref{eq:liley1})-(\ref{eq:liley4}) and (\ref{eq:liley6}) are studied in the present paper.
We note for the interested reader that a thorough 
discussion of realistic axonal velocity distributions in MFMs has been
recently published by Bojak and Liley \cite{ref:BL10} and that
Liley et al. \cite{ref:LCD02} contains more in-depth discussions of the physiological and 
mathematical arguments that lead to Eqs.~(\ref{eq:liley1})-(\ref{eq:liley5}). Also,
parameters used in the equations and their physiological ranges are shown in
Tab.~\ref{tab:allparms}.

\begin{sidewaystable}[htbf]
\begin{tabular}{lllll}

\hline\noalign{\smallskip}
Parameter & Definition & Minimum & Maximum & Units  \\
\hline
$h^r_k$ & resting membrane potential & $-80$ & $-60$ & mV\\
$\tau_k$ & passive membrane decay time & $5$ & $150$ & ms\\
$h^{\mathrm{eq}}_{ek}$ & excitatory reversal potential & $-20$ & $10$ & mV\\
$h^{\mathrm{eq}}_{ik}$ & inhibitory reversal potential & $-90$ & $ h^r_k-5$${}^\dagger$ & mV\\
$\Gamma_{ek}$ & EPSP peak amplitude & $0.1$ & $2.0$ & mV\\
$\Gamma_{ik}$ & IPSP peak amplitude & $0.1$ & $2.0$ & mV\\
$\gamma_{ek}$& EPSP characteristic rate constant$^\ddagger$ & $100$ & $1,000$ & $\mathrm{s}^{-1}$\\
$\gamma_{ik}$& IPSP characteristic rate constant$^\ddagger$ & $10$ & $500$ & $\mathrm{s}^{-1}$\\
$N^\alpha_{ee}$ & no.\ of cortico-cortical synapses, target excitatory & $2000$ & $5000$ & --\\
$N^\alpha_{ei}$ & no.\ of cortico-cortical synapses, target inhibitory & $1000$ & $3000$ & --\\
$N^\beta_{ek}$ & no.\ of excitatory intracortical synapses & $2000$ & $5000$ & --\\
$N^\beta_{ek}$ & no.\ of inhibitory intracortical synapses  & $100$ & $1000$ & --\\
$v$ & axonal conduction velocity & $100$ & $1,000$ & $\mathrm{cm}\,\mathrm{s}^{-1} $\\
$1/\Lambda$ & decay scale of cortico-cortical connectivity  & $1$ & $10$ & cm \\
$S^{\mathrm{max}}_k$ & maximum firing rate & $0.05$ & $0.5$  & $\mathrm{ms}^{-1}$ \\
$\mu_k$ & firing threshold &  $-55$ & $-40$ & mV \\
$\sigma_k$ & standard deviation of firing threshold  & $2$ & $7$ & mV \\
$p_{ek}$ & extracortical synaptic input rate & $0$ & $10,000$ & $\mathrm{s}^{-1}$ \\
\hline
\end{tabular}
\caption{{\bf Physiological ranges for spatially averaged model parameters.}
Common ranges does not imply shared values for different types
$k=e,i$ of neuronal target populations, hence physiological parameter space is 32-dimensional.
For bulk oscillations and for absent inhibitory extracortical input, $v\Lambda$ and 
$\Gamma_{ik}N^{\beta}_{ik}$ occur only in combination: this makes parameter space 
effectively 29-dimensional. ${}^\dagger$The upper limit
for $h^{eq}_{ie}$ should be $h_e^r-5~\mathrm{mV}$, but was erroneously set to
$h_i^r-5~\mathrm{mV}$ in Ref.~\cite{ref:BL05}. Since power spectra were obtained
around a chosen fixed point with $h_e^*> h^{eq}_{ie}$ for all parameter sets, this
has no significant consequences. $^\ddagger$Equivalently time to peak PSP amplitude given
by $1/\gamma_{lk}$.}
\label{tab:allparms}
\end{sidewaystable}

\section{Methods}

\subsection*{Parameter sets and bifurcation parameters}

A large number of parameter sets is 
considered in this study, representing a significant sample 
of the physiologically admissible parameter space of Liley's MFM. 
To generate such realistic parameter sets, the properties of
the model have been studied around fixed points 
of linearized version of Eqs.~(\ref{eq:liley1})-(\ref{eq:liley5}).  
Biologically plausible selection criteria were then
invoked concerning the properties of
the simulated EEG and the firing rates of 
their associated neuronal populations. 
Limits on the admissible ranges of parameter values were also applied, 
so that they are consistent with well-established neurophysiological 
features of mammalian cortex. Details on the methods
and the constraints applied in this search can be found in Bojak and Liley \cite{ref:BL05}.
As a result, a collection of 73,454 physiologically
meaningful parameter sets
is available and considered
in this study as representative of physiologically  
realistic parametric instantiations of model cortical activity in Liley's MFM.
Bifurcation theory is
an essential tool for discussing the qualitative changes that
solutions of dynamical systems undergo under variations of their parameters. 
Loosely speaking, this type of analysis provides a picture similar to a phase diagram
for physical systems, so that changes in the
properties of the solutions (i.e., bifurcations) may be considered as
a generalization of phase transitions \cite{ref:SS10}. This
allows us to explore parametric boundaries between qualitatively
different activities.

\begin{figure}[htbf]
\begin{center}
\includegraphics[width=\textwidth]{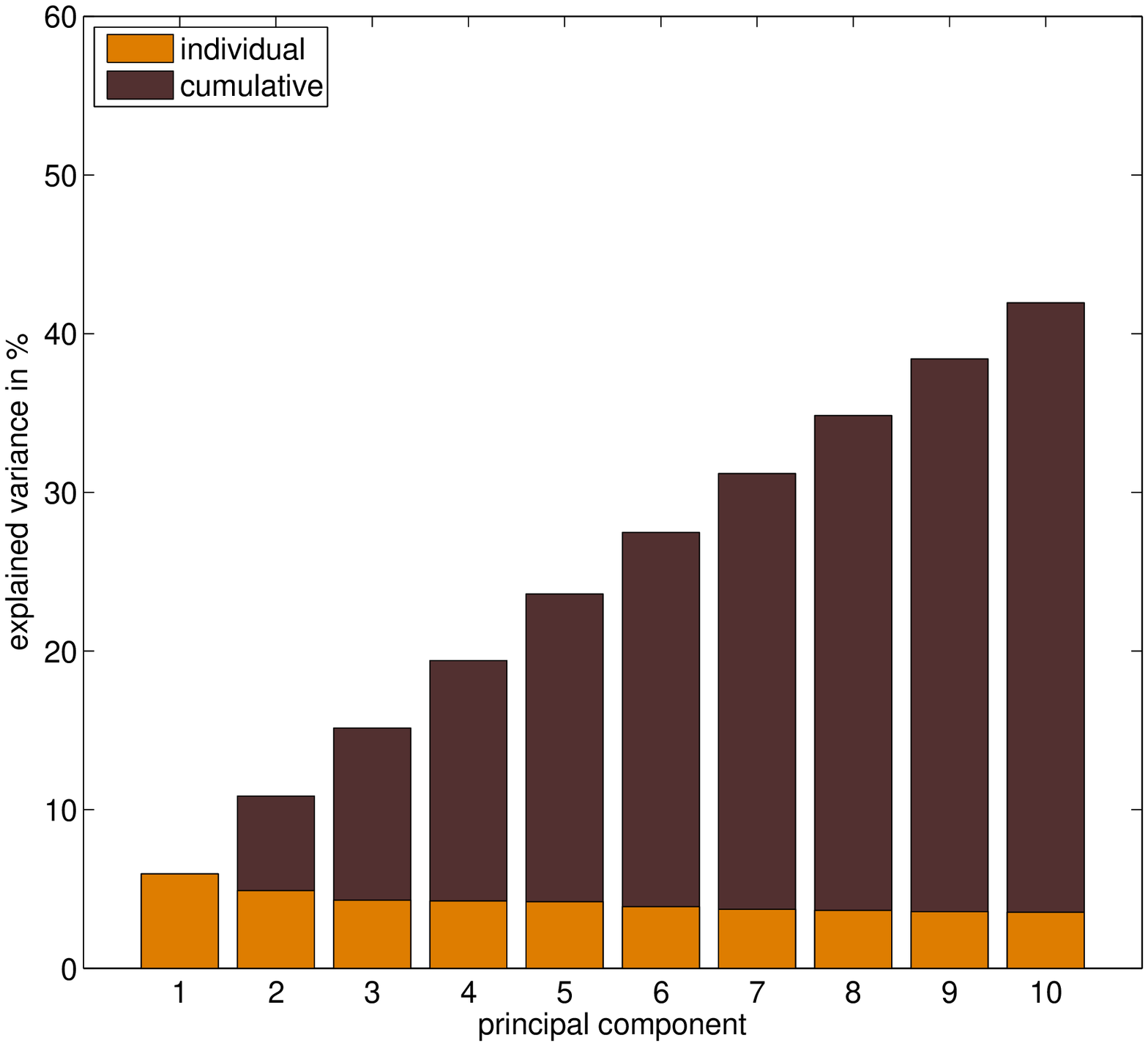}
\end{center}
\caption{{\bf Principal Component Analysis of the model parameter space.} 
Fraction of the total variance of parameters explained by the first ten principal components
for all 73,454 sets. Parameter space has 32 physiological dimensions, of which 29 are effectively
distinguished here, c.f.\ Tab.~\ref{tab:allparms}.
Orange bars indicate individual contributions, brown ones their cumulative 
sum. The evident absence of strong linear correlations within parameters means
parameter space cannot be reduced easily to a lower number of effective degrees of freedom.}
\label{fig:pareto}
\end{figure}

One of the a priori strongest motivations for the use of this approach
is represented by the failure of standard methods of statistical linear analysis,
which are unable to supply relevant information about the 
model's 29-dimensional parameter space.
An example of this is illustrated in Fig.~\ref{fig:pareto}. 
As a result of a Principal Component Analysis (PCA)
over the parameters 
of the 73,454 sets, the percentage of total variance explained   
by the first ten PCs is plotted, individually and cumulatively.
Essentially, if the parameter space presented strongly linear correlations
among the variables, it would be possible to identify a limited number of
components capturing the vast majority of the variance, indicating that the significant
degrees of freedom in the space were in fact much less than its dimension. This
is clearly not the case, since the first ten components cumulatively account for only about 42\% of
the total variance, with low and very similar singular contributions. It is not possible
to find a linear combination of parameters that expresses the ``bulk''
properties of the parameter space. If we want to extract valuable statistical information,
a different, non standard method of analysis is needed. 
On this basis, the aim of this paper is to highlight 
correlates between dynamical ``snapshots'' of the model,
i.e. two parameters bifurcation plots, and spectral and parametric properties of Liley's MFM.

In choosing the parameters to vary for our bifurcation analysis,
we are guided by the existence of literature 
pointing to inhibition being a sensitive locus of control
of brain dynamical responses. For instance, it is well-known that
pharmacological agents alter cortical activity primarily through 
acting on inhibition \cite{ref:RA04}. Factors that modify inhibitory PSP amplitudes and anatomical coupling
strengths of the inhibitory neural population are thus selected as continuation parameters. 
This is in line with previous work \cite{ref:BL07} and implies changes in Equation \ref{eq:liley4}
for the inhibitory PSP amplitudes
\begin{equation}
\Gamma_{ie} \rightarrow R\Gamma_{ie}\;,\qquad\Gamma_{ii} \rightarrow R\Gamma_{ii}\;,
\end{equation}
and local inhibitory-inhibitory connectivity 
\begin{equation}
N^{\beta}_{ii} \rightarrow k N^{\beta}_{ii}\;,
\end{equation}
where
$R$ and $k$ are chosen as our bifurcation parameters.

This choice is not exclusive and
other parameters or factors could  have reasonably been varied,  
possibly highlighting other poignant features that have not emerged
through continuations in $R$ and $k$. Nonetheless, these two parameters
allow biologically plausible control:
$R$ changes the amplitude of the inhibitory PSPs for both
the excitatory and inhibitory populations,
which for example general anesthetics and sedative agents 
can affect to a large degree. 
At the same time inhibitory-inhibitory connection strengths are altered via $k$, 
which changes the impact on the inhibitory feedback loop largely responsible
for generating specific oscillations \cite{ref:LCD02}. Note that
since $p_{ik}\equiv 0$ in Eq.~\ref{eq:liley4}, the effect of
inhibition on the excitatory population scales precisely with $R$,
but on the inhibitory population with $R\cdot k$. Thus, for example,
setting $k=1/R$ would simulate an exclusive
impact on the excitatory population. Essentially then,
we can control overall inhibitory ``potency'' with $R$ and ``specificity'' with $k$.

\subsection*{One and two parameters bifurcation analysis procedure}

There are two types of bifurcation plots that will be discussed:
one parameter (\textsf{1par}) plots in $R$ and two parameters (\textsf{2par})
graphs in $R$ and $k$. 
Let us first explain how the diagrams are obtained: 
we integrate Eqs.~(\ref{eq:liley1})-(\ref{eq:liley4}) and (\ref{eq:liley6})
for a chosen parameter set
up to an equilibrium solution for $R = 1$ and $k = 1$, i.e., in the absence of changes to inhibition. 
Initial conditions are given at $t=0$ by $h_e(0)
= h_i(0) = - 70~\mathrm{mV}$, while all the other state variables are set to zero. Once an 
equilibrium solution is found, it is continued in $R$ and a \textsf{1par} 
plot is obtained. In all the $R\/$-plots there are
always two types of codimension one (\textsf{cod1}) bifurcations present:
saddle-node and Hopf points, respectively.
The presence of Hopf bifurcations gives us some information about the character of
sustained oscillations produced by the action of inhibition, which we will discuss below.

The saddle-node and Hopf bifurcations are subsequently continued in parameters $R$ and $k$.
In our diagrams, lines describing the loci of saddle-nodes are labeled \textsf{sn}, and similarly
for the Hopf points \textsf{hb}. 

Physiological limits for minimal and maximal parameters intervals are set at 
$R_{min} = k_{min}= 0.75$, $R_{max} = 2.00$ and $k_{max} = 1.25$.  
Higher or lower values of $R, k$ are considered to
alter inhibitory PSP amplitudes and
inhibitory-inhibitory connectivity strengths either to
unphysiological values or to unusually large or small
figures, which, for example, do not capture the commonly observed
effects of pharmacological agents on the cortex. However, as discussed below,
continuations are performed till particular features of the bifurcation
topology become apparent. Typically, this occurs far outside of the
physiological range.
Bifurcation analysis is performed here with the widely used
continuation software packages AUTO-07P and MATCONT \cite{ref:Doe81,ref:DGK03}.

\section{Results}

\subsection*{Two parameters bifurcation patterns}

To establish our findings, we proceeded with a two-step analysis. 
First, \textsf{2par} bifurcation plots for a reduced sample of sets were obtained and
analyzed, visually and numerically. Analysis of plots in $R, k$ for 405 randomly selected sets out
of the 73,454 in the generated sample
revealed the existence of two unique, and clearly distinct, 
categories of patterns which recurrently appear. Our primary aim in this subsection
is to describe these archetypal patterns 
and discuss their shared and distinct features. In the next subsection, 
a third parameter ($p_{ei}$) is varied and this shows that the two families are
``transformable'': it is possible to metamorphose one family into the other,
by varying exogenous stimuli. Once the two families
are characterized, generalization through the whole batch is straightforward
and global features of the physiological parameter space naturally emerge.

\begin{figure}[htbf]
\begin{center}
\includegraphics[width=\textwidth]{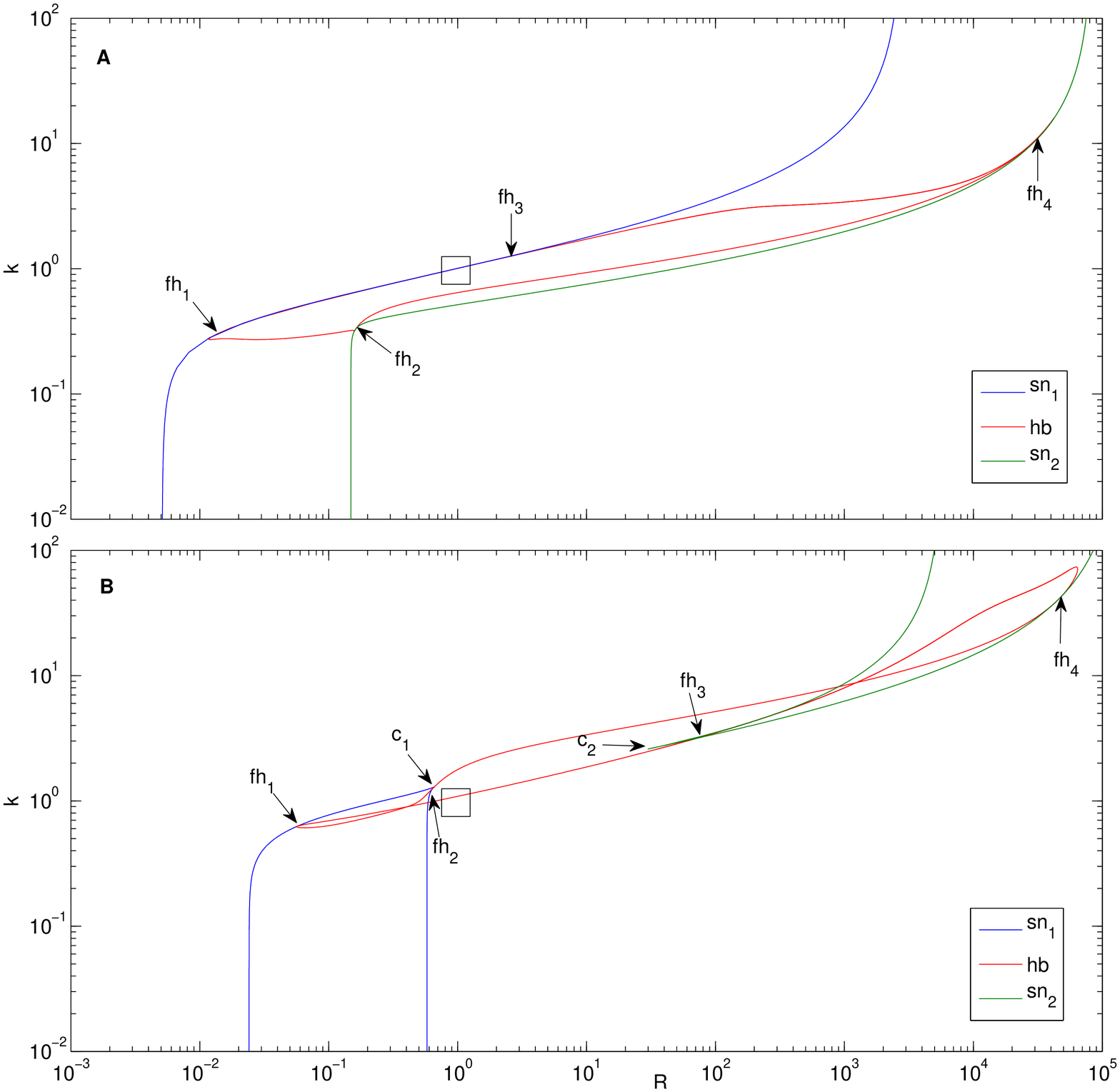}
\end{center}
\caption{{\bf Typical bifurcation structure of families.} Continuations
for two distinct parameter sets belonging to F1 ({\bf A}) and F2 ({\bf B}),
respectively.
Note in {\bf B} that fold-Hopf points $\mathsf{fh}_{2}$ and $\mathsf{fh}_{3}$ of F2 occur close to 
the cusps $\mathsf{c}_{2}$ and $\mathsf{c}_{1}$, respectively.
The two saddle-node lines are shown in blue $\mathsf{sn}_{1}$
and green $\mathsf{sn}_{2}$, respectively, and the Hopf line \textsf{hb} in red.
The number of equilibria is $3$
for areas bounded by the saddle-node lines \textsf{sn}, and $1$ outside.
The rectangle shows the physiologically relevant range of $Rk\/$-space.
}
\label{fig:fams}
\end{figure}

The general character and topology of the \textsf{2par} plots
for Liley's MFM is reproduced in Fig.~\ref{fig:fams}.
These represent two archetypal $Rk$-plots for so-labelled Family~1 (F1) and
Family~2 (F2) types of sets. The major difference among the two families is 
represented by how the \textsf{sn} lines 
are organized and divide the parameter space in areas with different
dynamical properties. In general, because of the form of the MFM's equations, 
a region bound by \textsf{sn} lines has always three equilibria,
whereas the region outside has only one. For F1, see Fig.~\ref{fig:fams}A,
\textsf{sn} lines travel almost parallel for the whole range of the space, 
and do not form cusp points (\textsf{c}), dividing the 
space essentially in three major regions: one inside \textsf{sn} lines 
and two outside. In contrast, F2, see Fig.~\ref{fig:fams}B, is characterized by the presence of two 
cusp points $\mathsf{c}_{1}$ and $\mathsf{c}_{2}$, so that the region containing three equilibria is the union of two separated
``wedge-shaped'' areas, with cusps as their vertices. 
The way \textsf{hb} and \textsf{sn} lines interact is, for the majority of cases,
by so-called fold-Hopf (\textsf{fh}) points, which appear as tangencies between 
such lines. 
Other cases we saw but are
not present in Fig.~\ref{fig:fams} happen
when \textsf{hb} lines connect two 
Bogdanov-Takens (\textsf{bt}) points. These instead
correspond in the diagram to points where a \textsf{hb} line 
terminates on a \textsf{sn} line.
In general, these \textsf{bt} points 
occur on the left, lower wedge of the \textsf{sn} lines, close
to \textsf{fh}$_{1}$ for F1 and F2. We will give examples of
these bifurcations in the next subsection.

For both families, there are also other cases that slightly differ from Figure \ref{fig:fams},
although the structure of the \textsf{sn} lines never changes. 
Variations are minimal and local, i.e. one of the \textsf{fh} points can be missing
or a small difference in the shape of the \textsf{hb} lines can be 
present, but the essential structure 
of the plot is not altered. What instead shows high variation, 
independently of the family, is the extension of the plots, i.e.,
the values of $R$ and $k$ at which 
bifurcation points appear. 
It is important to realize that the physiological region
is typically a tiny area of the plots, as indicated by the rectangles in Fig.~\ref{fig:fams}.
There are plots, as those shown in Fig.~\ref{fig:fams}, 
where the whole pattern is contained within $3 - 4$ orders of magnitude
of the bifurcation parameters, but there are also a large number of patterns for 
which \textsf{hb} and \textsf{sn} lines extend up to immense values, like $R \approx 10^9$. 
This is an interesting, important point: scales vary widely, but the nature 
and structure of the patterns are invariant, and this 
is appreciated only if plots are considered in their 
complete extension in $Rk$-space. Moreover, it appears that
similar bifurcations maintain the same type of normal forms among
equivalent patterns, no matter what scale is involved. 
This implies that the dynamics in the proximity of bifurcations
is topologically equivalent, i.e. qualitatively very similar, among plots in the same family. 
Thus, parameter sets in the same family share not only the topology of the saddle-node
bifurcation lines, but also yield qualitatively similar dynamics close to the codimension
two points.

\subsection*{Thalamic input and pattern family metamorphosis}

The magnitude of extracortical (thalamic) input to inhibitory neurons $p_{ei}$
is chosen as a third continuation parameter for investigating 
the transition between families. Our choice of
$p_{ei}$ has an intuitive meaning: unlike $\Gamma_{ie}$ and $\Gamma_{ii}$,
whose variation represent endogenously driven effects on the cortex, changes in the thalamic input
capture the action of exogenous causes. Furthermore, $p_{ei}$ excites the 
inhibitory population and hence can be expected to affect significantly the dynamical repertoire
of families based on variations of inhibitory PSPs. 

\begin{figure}[htbf]
\begin{center}
\includegraphics[width=\linewidth]{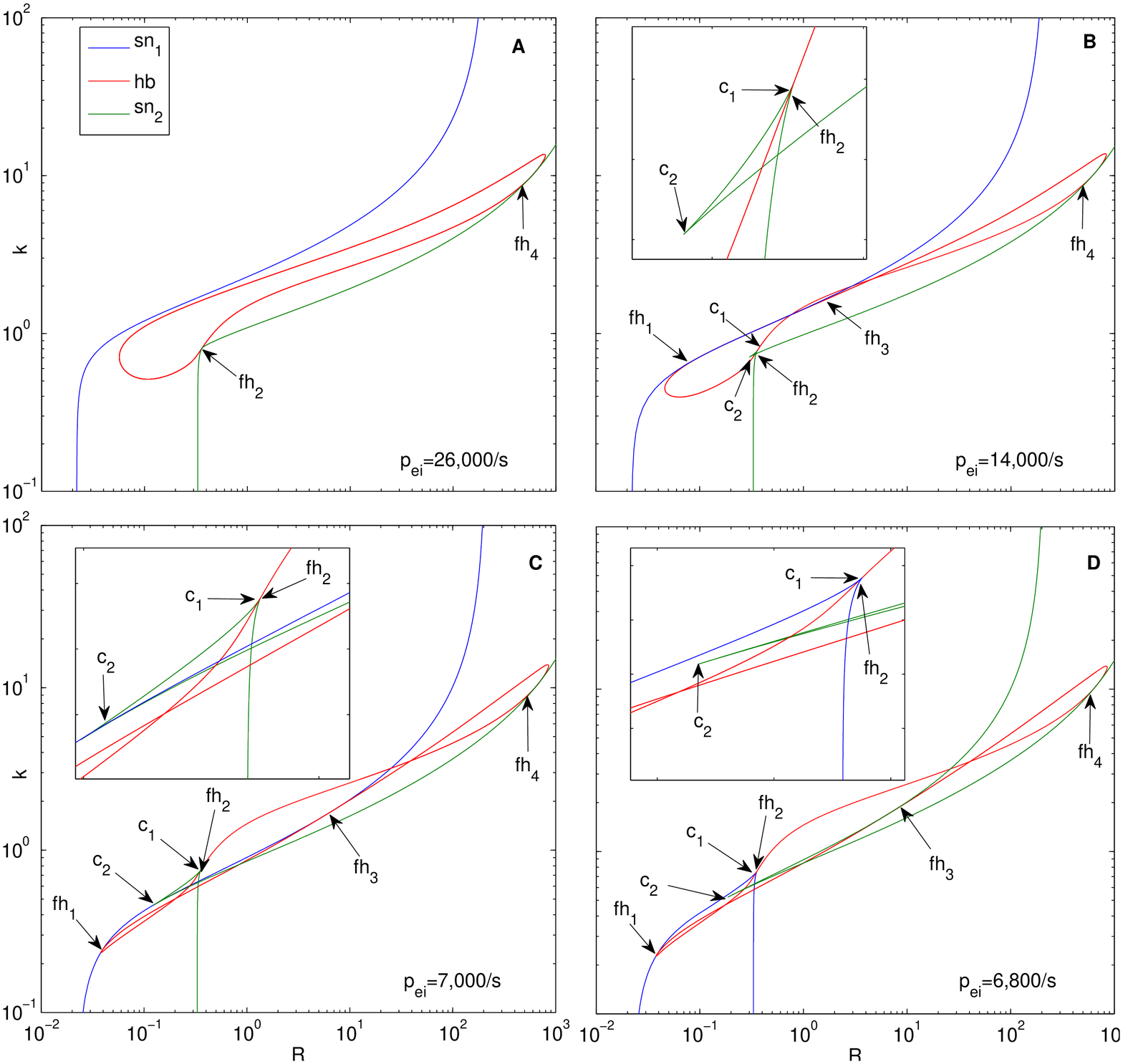}
\end{center}
\caption{{\bf Exogenous effects: inter-familiar transition.} Metamorphosis from F1 to 
F2 through two swallow tail bifurcations, induced by decreasing thalamic excitation
to inhibitory neurons:
{\bf A} F1 at $p_{ei} = 26,000/\mathrm{s}$,
{\bf B} F1 at $p_{ei} = 14,000/\mathrm{s}$,
{\bf C} just before bifurcation at $p_{ei} = 7,000/\mathrm{s}$,
and {\bf D} just after bifurcation $p_{ei} = 6,800/\mathrm{s}$. 
The separation between the cusps $\mathsf{c}_{1}$ and $\mathsf{c}_{2}$ further increases
for even lower values of $p_{ei}$, and a pattern identical
to the typical F2 in Fig.~\ref{fig:fams}B appears at $p_{ei} \approx 6,000/\mathrm{s}$.
Analogous transitions occur for increasing thalamic excitation to excitatory neurons ($p_{ee}$).
The overwhelming majority of the inspected sample of 405 sets metamorphose similarly.}
\label{fig:F1toF2}
\end{figure}

Let us discuss a prototypical transition from F1 to F2. 
To follow our description the reader will need to refer to  Fig.~\ref{fig:F1toF2} . 
We start at a very large and unphysiological value 
of $p_{ei} = 26,000/\mathrm{s}$, indicated in 
Fig.~\ref{fig:F1toF2}A, where a representative of F1 is present. Compared to the
(standard) diagram depicted in Fig.~\ref{fig:fams}, two \textsf{fh} points
on the left wedge of line $\mathsf{sn}_1$ are not yet formed. 
When the input is decreased to $p_{ei} = 14,000/\mathrm{s}$ in Fig.~\ref{fig:F1toF2}B,
thus approaching the physiologically plausible range of 
$0/\mathrm{s}\leq p_{ei}\leq 10,000/\mathrm{s}$,
the \textsf{hb} branch extends towards lower values of $R$,
creating the $\mathsf{fh}_{1}$ point, as found in the
standard F1 type in Fig.~\ref{fig:fams}. Analogously, 
the $\mathsf{fh}_{3}$ point is created at larger values 
of $R$ and $k$. At the same time, on the right wedge 
$\mathsf{sn}_2$, two cusps $\mathsf{c}_{1}$ and $\mathsf{c}_{2}$ appear
in Fig.~\ref{fig:F1toF2}B. This happens in a so-called {\em swallow tail}
bifurcation (see, e.g., Ref.~\cite{ref:arnold}). All equilibria belong to 
the same branch, meaning that they are connected
at the saddle-node bifurcations, at which the branch is folded. 
At this swallow tail bifurcation, a pair of saddle-node bifurcations 
appears on the unique branch of equilibria.
Notice in the inset that the $\textsf{fh}_{2}$ point 
is now on the right wedge of the little triangle formed by the two cusps
$\mathsf{c}_{1}$ and $\mathsf{c}_{2}$. Moreover, the single branch \textsf{hb} has lost
its previous elliptical shape, and has produced intersections.
However, these are not bifurcation points since
only one unique branch of \textsf{hb} lines is present, and the diagram
still represents a type of plot belonging to F1. 
The emergence of the two cusps is not considered a major change, 
since the topological characteristic that
we choose in discriminating F1 from F2 is the global shape of the 
\textsf{sn} lines. When separate $\textsf{sn}_1$ and $\textsf{sn}_2$ do not
intersect but keep almost parallel, as in Fig.~\ref{fig:F1toF2}A and B,
a diagram of F1 type is still present.

As the extracortical input is further decreased to $p_{ei} = 7,000/\mathrm{s}$
in Fig.~\ref{fig:F1toF2}C, the triangle of line 
$\mathsf{sn}_2$ with the two cusps $\mathsf{c}_{1}$ and $\mathsf{c}_{2}$ expands
and moves to the left of the diagram,
intersecting the opposite $\mathsf{sn}_1$ branch (see also the inset in 
Fig.~\ref{fig:F1toF2}C). Now $\mathsf{c}_{2}$ is 
overlapping the line of saddle-nodes $\mathsf{sn}_1$, and the cusp $\mathsf{c}_{1}$ 
has crossed $\mathsf{sn}_1$ and is above it. Notice, in particular,
the vicinity of the blue and green lines of \textsf{sn} in the inset: 
this represents a stage very close to the transition of F1 into F2. In fact, 
at smaller $p_{ei}$ values, the cusp $\mathsf{c}_{2}$ further lowers and
separation between the two branches of \textsf{sn} further diminishes, so that
increasing portions of $\mathsf{sn}_1$ and $\mathsf{sn}_2$ become 
nearly tangent. 
This condition persists up to a value of $p_{ei}$ where \textsf{sn} lines finally
exchange, so that at $p_{ei} = 6,800/\mathrm{s}$ the 
two cusps no longer belong to the same branch,
as Fig.~\ref{fig:F1toF2}D and its inset show. The cusp $\mathsf{c}_{1}$ is now part of $\mathsf{sn}_1$ and 
$\mathsf{c}_{2}$ belongs to $\mathsf{sn}_2$: a second {\em swallow tail} bifurcation 
has taken place, signaling the appearance of an F2 type of diagram. In other words, 
two regions of the $Rk$-space having triangular 
shapes have formed and overlapped, and $\mathsf{c}_{2}$ now tends 
to move towards the right as $p_{ei}$ is decreased, 
separating further from $\mathsf{c}_{1}$. A further decrease of the 
thalamic input will eventually result in the same plot as in Fig.~\ref{fig:fams}B, 
i.e., the prototypical F2 with two disjoint V and $\Lambda$-shaped curves of saddle-nodes.

A physiologically relevant observation can be made at this point.
One can investigate the effect of varying the excitatory thalamic input to either
inhibitory ($p_{ei}$) or excitatory ($p_{ee}$) cortical neurons.
It turns out that there is a relatively robust family-specific effect of thalamic excitation,
which is opposed for inhibitory and excitatory cortical targets:
F1 (F2) is the preferred family at high (low) values of $p_{ei}$, whereas
F1 (F2) is the preferred family at low (high) values of $p_{ee}$. The incidence of family metamorphosis
for variations {\em within} the physiological interval of thalamic input 
$0/\mathrm{s}~-~10,000/\mathrm{s}$ is quite strong: about 56\% of the 73,454
sets undergo a transition when $p_{ei}$ is varied and approximately 42\% when $p_{ee}$ is. Once again,
inhibition seems to be somewhat more effective in controlling the global dynamical properties of
the cortex, in this case as driven by excitatory exogenous inputs.

\subsection*{Changes within families due to thalamic input}

To complete our analysis of \textsf{2par} plots,
let us make some remarks on the other instantiations
of F1 and F2 that have not been illustrated yet, but do occur in the 405 analyzed sets.
These diagrams are only slight, local variations of the ones shown so far,
and still exhibit the global common characteristics typical of the respective
families. As an example, patterns belonging to F1 very similar to 
Fig.~\ref{fig:F1toF2}A are present, 
but they miss the $\mathsf{fh}_{3}$ point. This type
of diagram is still considered as
part of F1, because the structure of its \textsf{sn} lines does not change.
If we start looking at the transition between F1 and F2 from this particular plot,
rather than from Fig.~\ref{fig:F1toF2}A,
the fate of the pattern as $p_{ei}$ is decreased is the same, i.e.,
it metamorphoses into an F2 diagram.

\begin{figure}[htbf]
\begin{center}
\includegraphics[width=\linewidth]{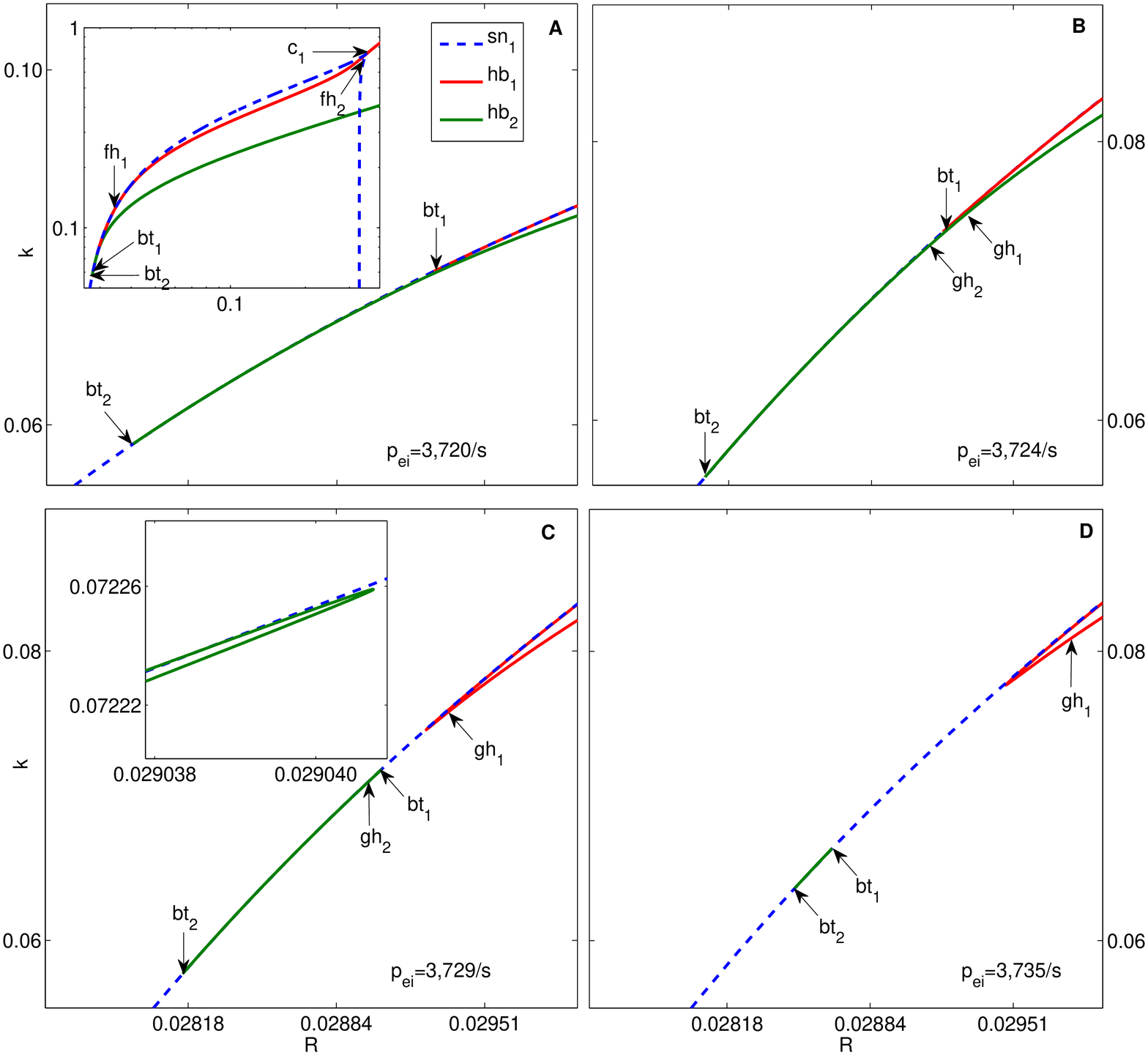}
\end{center}
\caption{{\bf Exogenous effects: annihilation of \textsf{bt} points.}
Increasing thalamic excitation to inhibitory neurons results in the annihilation of the \textsf{bt} points:
{\bf A} $p_{ei} = 3,720/\mathrm{s}$,
{\bf B} $p_{ei} = 3,724/\mathrm{s}$, {\bf C} $p_{ei} = 3,729/\mathrm{s}$ and {\bf D} $p_{ei} = 3,735/\mathrm{s}$. 
The inset in {\bf A} shows the triangular region of F2 where the bifurcations take place. 
Further increase of $p_{ei}$ shrinks \textsf{hb}$_{1}$ up to the coalescence and disappearance of the \textsf{bt} points.
Note that the transformations induced by the thalamic input
occur over a very small range of $p_{ei}$ and within even smaller intervals for $R$ and $k$,
demonstrating the sensitivity of the parameters. To improve visualization,
continuation lines have a larger linewidth and the \textsf{sn} branch is dashed.
}
\label{fig:AnnBTs1}
\end{figure}

There is one final example of possible variations, 
regarding the appearance of Bogdanov-Takens (\textsf{bt}) points,
which we anticipated in the previous subsection.
Remarkably, these changes can be accessed as well
by varying the same parameter $p_{ei}$: exogenous
agents do not only produce transitions between families, 
but can also trigger modifications {\it within} the same family. 
We study here the annihilation of a couple of \textsf{bt} points, 
by increasing $p_{ei}$ for the same parameter set employed in
Fig.~\ref{fig:F1toF2}. Once again the reader is asked
to look closely at Fig.~\ref{fig:AnnBTs1} and its insets to follow our description. 
If we start at a low, physiological $p_{ei} = 3,720/\mathrm{s}$, 
two Bogdanov-Takens points, $\mathsf{bt}_1$ and $\mathsf{bt}_2$, are present just below 
$\mathsf{fh}_{1}$ and we want to demonstrate how they disappear for larger
$p_{ei}$.
These \textsf{bt}s occur on the left wedge of the \textsf{sn} line culminating in the
cusp point $\mathsf{c}_{1}$. Notice also that the branch $\mathsf{hb}_{1}$ 
is responsible for the point $\mathsf{fh}_{2}$ close to $\mathsf{c}_{1}$, 
in analogy with the F2 in Fig.~\ref{fig:fams}. The similarity between the two plots
is evident. As the thalamic input is increased to $p_{ei} = 3,724/\mathrm{s}$ in 
Fig.~\ref{fig:AnnBTs1}B,
a couple of generalized Hopf (\textsf{gh})
points appear on $\mathsf{hb}_{1}$. These points
signal a change in the \textsf{hb} line
from subcritical to supercritical, i.e. the Hopf bifurcation
associated to the branch switches from a hard to a soft type \cite{ref:Kuz04}. 

If $p_{ei}$ is made even larger, the two \textsf{hb} lines come closer and closer,
up to a value of the input where they glue together and separate, 
as shown in Fig.~\ref{fig:AnnBTs1}C for $p_{ei} = 3,729/\mathrm{s}$. After this, $\mathsf{hb}_{1}$ 
acquires a large portion of what previously was $\mathsf{hb}_{2}$, with the
latter now extending over a much smaller distance. 
The split occurs in between $\mathsf{gh}_{1}$ and $\mathsf{gh}_{2}$, each of which 
now belongs to a separate Hopf branch: the tiny $\mathsf{hb}_{2}$ reaches
from $\mathsf{bt}_1$ to $\mathsf{bt}_2$ and the long $\mathsf{hb}_{1}$ bends in proximity of $\mathsf{gh}_{2}$ 
and continues to larger $R$ and $k$ values. Notice in the inset of
Fig.~\ref{fig:AnnBTs1}C 
that also $\mathsf{hb}_{2}$ (in green) bends close to $\mathsf{gh}_{1}$ and heads back to $\mathsf{bt}_1$,
where it stops. Then, when $p_{ei}$ is further incremented,
see Fig.~\ref{fig:AnnBTs1}D with $p_{ei} = 3,735/\mathrm{s}$,
$\mathsf{gh}_{2}$ annihilates and the smaller branch $\mathsf{hb}_{2}$ 
further shrinks, with the \textsf{bt} points getting closer and closer
as $p_{ei}$ grows. In the end, these two \textsf{bt} points coalesce, 
and we are left with only a single branch, 
as in Fig.~\ref{fig:fams}B. We note that this line of \textsf{hb}
delimited by $\mathsf{bt}_1$ and $\mathsf{bt}_2$ only exists at
positive values of $R$ and $k$. For the interested reader, we point out
that the plot in Fig.~\ref{fig:AnnBTs1}A is responsible
for chaos in an unphysiological range of $R$ and $k$, 
and its route through a so-called
homoclinic doubling cascade is analyzed in detail 
in \cite{ref:FvVB+08}.

Let us conclude with some observations. First, the tiny $\mathsf{hb}_{2}$ branch 
exists for a very limited interval of $p_{ei}$ and 
this is reflected in its very low occurrence in the 405 sets: it
was found only twice. Second, an identical mechanism for F1 is also present,
where the creation or annihilation of \textsf{bt} points takes place on the 
left $\mathsf{sn}_1$ lines in Fig.~\ref{fig:fams}A. 
Finally, the diagrams
that have been presented for the annihilation of \textsf{bt} points
in F2, the similar 
plots for F1 (not shown), and the plots in Figs.~\ref{fig:fams} and \ref{fig:F1toF2}
with their small, local variations, are the only ones
that have emerged in our analysis of the reduced batch
of 405 sets. Within the previously discussed assumptions, 
we conclude that these diagrams constitute all the possible modeled cortical responses
to inhibition that are physiologically meaningful in Liley's MFM.

\subsection*{Partition of parameter space: families reaction to anesthetics and distributions of parameters}

Having shown the properties of F1 and F2 types of diagrams, and concluded that any of the 405 plots of the 
preliminary sets falls into one of the two, we process
the whole 73,454 sets in search of global, statistically relevant, correlations between family membership and
parameter distributions. 
The approximation we introduce to cope with the large number of sets in an automated fashion  
is simple: F1 and F2 are those sets respectively
without and with separate lines of \textsf{sn}, which are divided
by the exchange of \textsf{sn} branches illustrated in Figs.~\ref{fig:F1toF2}C-D. 
Hence, Figs.~\ref{fig:fams}A and \ref{fig:F1toF2}A-C
are considered as F1, Figs.~\ref{fig:fams}B and \ref{fig:F1toF2}D as F2. 
The second swallow tail bifurcation mentioned previously
is the boundary between the so-defined F1 and F2.
The algorithm we employed searches for the position and number of cusps on \textsf{sn} lines,
and assigns the family type consequently. Equilibration, continuation and attribution of a single set to its 
family takes on average $9$ seconds on an everyday desktop computer.

The first, global difference between families we wish to show lies in their response to 
modeled anesthetic action. According to the methods explained in Bojak and Liley \cite{ref:BL05},
it is possible to simulate the increase or decrease in EEG power when anesthetics are induced,
and extract the ratio between the power of the anesthetized cortex and the cortex at rest.
Herein power is predicted for the system at a stable equilibrium subject to fluctuations
induced by noise added to the thalamocortical input $p_{ee}$.
In line with clinical practice, they
adopted the minimum alveolar concentration (MAC) 
of anesthetic agent at 1 atm pressure as the reference
measure for describing the behavior of the anesthetized cortex.  
One MAC is the inspired anesthetic concentration needed to prevent movement
in 50\% of people to a noxious (surgical) stimulus.
For example, assuming that isoflurane, a common anesthetics, is employed during surgery, a patient would be 
maintained usually at $0.9$ to $2.2~\mathrm{MAC}$ isoflurane
in an oxygen-$70\%$ nitrous oxide mixture, or at $1.3-3~\mathrm{MAC}$
without the nitrous oxide. According to previous studies, it is assumed that
1 MAC is equal to a $1.17~\mathrm{vol\%}$ for isoflurane, which 
corresponds to $c = 0.243~\mathrm{mM}$ aqueous concentration \cite{ref:Map76,ref:FL96}. 

\begin{figure}[htbf]
\begin{center}
\includegraphics[width=\linewidth]{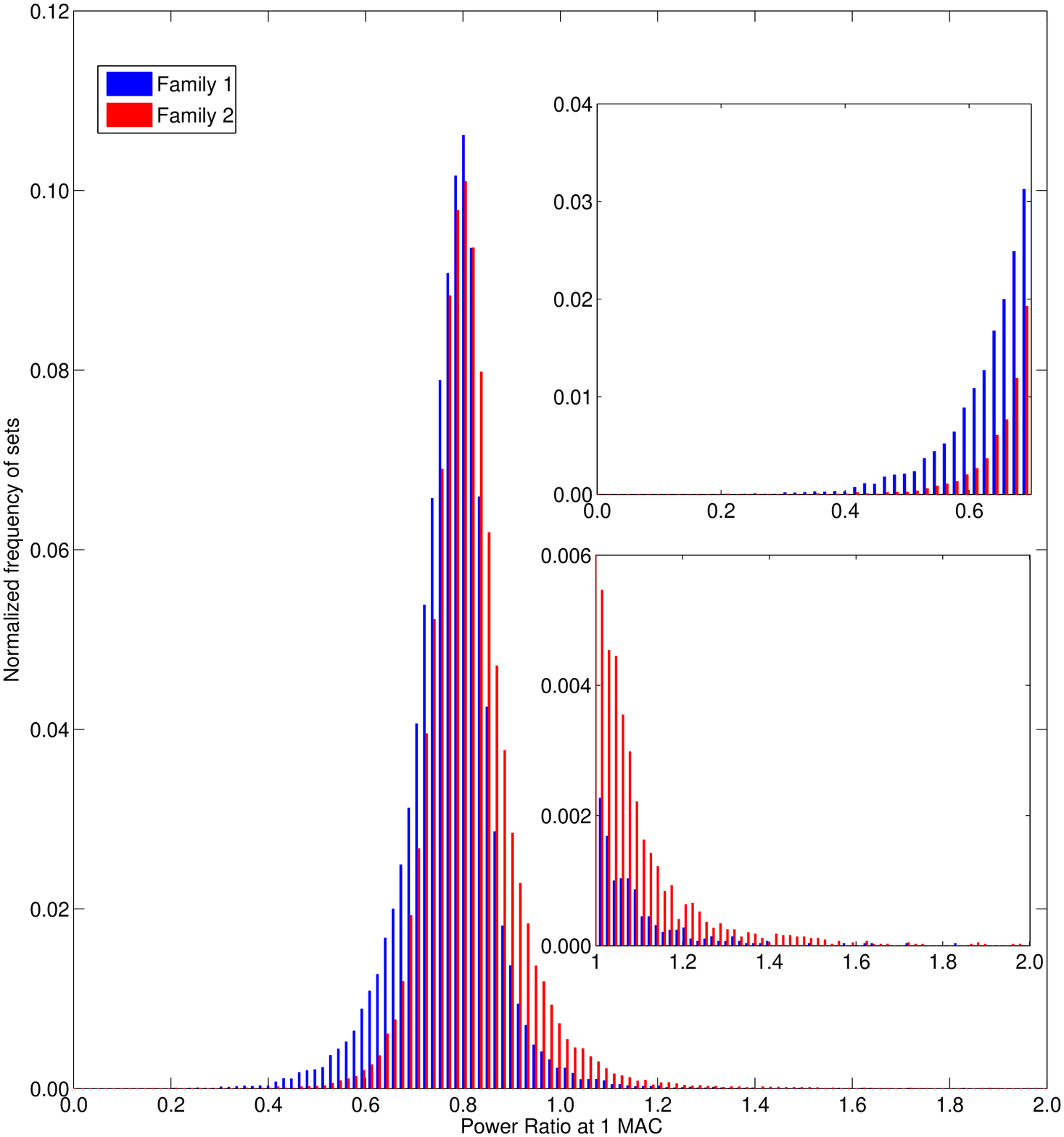}
\end{center}
\caption{{\bf Distributions of power ratios sorted by families.} 
The normalized distributions of power rations for the analyzed 73,454 sets, with F1 in blue and F2 in red, and magnifications of their tails in the insets.
There are 29,131 sets of F1 type (39.7\% of the total) and 44,323 sets of F2 type (60.3\%).
The total number of bins is 250 for power ratios $0-4$, of which only $0-2$ is shown.
Average and standard deviation of the power ratios
are $0.77 \pm 0.09$ for F1 and $0.82 \pm 0.09$ for F2.
Cumulative probabilities for the distributions at the tails are
$P($power ratio $\leq 0.7 | F1) = 0.161$,
$P($power ratio $\leq 0.7 | F2) = 0.060$,
and $P($power ratio $\geq 1.0 | F1) = 0.012$, $P($power ratio $\geq 1.0 | F2) = 0.037$.
}
\label{fig:allrel}
\end{figure}

The distributions of family responses to anesthetics according to power 
ratios are depicted in Fig.~\ref{fig:allrel}, 
with F2 types being in the majority in the whole batch,
since they turn out to be about one and half times as many as F1. 
The histograms look similar in shape but are shifted with respect to
one another. They are characterized by similar averages,
with a relative difference of only about $5\%$, but are strikingly 
different at the respective tails. At low ratios, F1 sets have a tendency of appearing
more frequently than F2 sets, and vice versa for high ratios.
For example, for power ratio values smaller than $0.7$, 
the cumulative probability for F1 is almost three times that of F2,
whereas for ratios bigger than $1.0$ the situation is reversed, with F2 having a 
cumulative probability more than three times that of F1. This 
shows that global, dynamical patterns have an interesting degree of
correlation with total EEG power when responding to modeled isoflurane,
whose mode of action is representative of most anesthetics.

Intuitively, anesthetized cortices 
should show a decrease in total EEG power when 
anesthetics are induced, as behavior is impaired as the cortex ``slows down'', with
power ratios at 1 MAC expected to be smaller than one. This is not always true, and it is not 
uncommon to rather observe a transient surge in EEG power that is clinically known 
as the biphasic response \cite{ref:BGC+97}. 
For the generated 73,454 sets, a power increase larger than 1.4 at 1~MAC, 
has been taken as indicative of a biphasic response in \cite{ref:BL05}, i.e., 
the simulated total power integrated over all frequencies
was at least 1.4 times larger with than without the presence of 1~MAC isoflurane.
Only 86 sets out of the whole batch show a biphasic response, of which 
around 90\% are of F2 type. This is interesting, but the statistics are too poor 
to draw a general conclusion about the relation of F2 to
the appearance of a biphasic response. It is possible that the predominance 
of F2 for high ratios is related to the cortical mechanism eliciting
a surge in power when anesthesia is present.

Besides responses to anesthetics, our family partition also
shows that some frequency distributions of parameters
exhibit significant differences for F1 and F2. In order to objectively assess
these differences, we use the following procedure: first, we construct
parameter frequency histograms with $n=20$ bins within the limits provided
for each parameter by Tab.~\ref{tab:allparms}, separately for the two families.
Two particular cases arise in this procedure.
On one hand for $h_{ie}^{eq}$ and $h_{ii}^{eq}$, which have flexible upper limits,
we choose the maximum possible one ($-65~\mathrm{mV}$) as bin limit.
On the other hand, $\gamma_{ii}$ occupies almost exclusively the lower part of
the allowed parameter range. For the sake of better accuracy with the same number of bins,
we set $120/\mathrm{s}$ as the upper bin limit for $\gamma_{ii}$. Parameter sets with larger values,
i.e. 12 out of 29,131 for F1 and 24 out of 44,323 for F2, are counted in the
highest bin. Second, we now compute the square root of the information radius
(also known as the Jensen-Shannon divergence) between the histograms
for each parameter $j$ -- $g^j_i$  for F1 and $h^j_i$  for F2 -- with $i=1,\ldots,n$ bins 
\begin{equation}
\begin{gathered}
d_{\mathrm{IR}}^j=\sqrt{\sum_{i=1}^n\left[\frac{1}{2}
\left(g^j_i\log_2 g^j_i+h^j_i\log_2 h^j_i\right)-a_i\log_2 a_i\right]}\;,\\
\mathrm{with}\quad a_i=\frac{1}{2}\left(g^j_i+h^j_i\right)\;,
\end{gathered}
\end{equation}
which is a proper metric for the similarity of discrete probability distributions
\cite{ref:ES03}. 
Note that for empty bins, one defines $0\log_2 0\equiv 0$. Identical
distributions have $d_{\mathrm{IR}}=0$, whereas maximal dissimilarity is indicated
by $d_{\mathrm{IR}}=1$. An example for maximal dissimilarity would be histograms $g^j_i=\delta_{ik}$
and $h^j_i=\delta_{il}$ with Kronecker deltas and $k\neq l$.

\begin{figure}[htbf]
\begin{center}
\includegraphics[width=\linewidth]{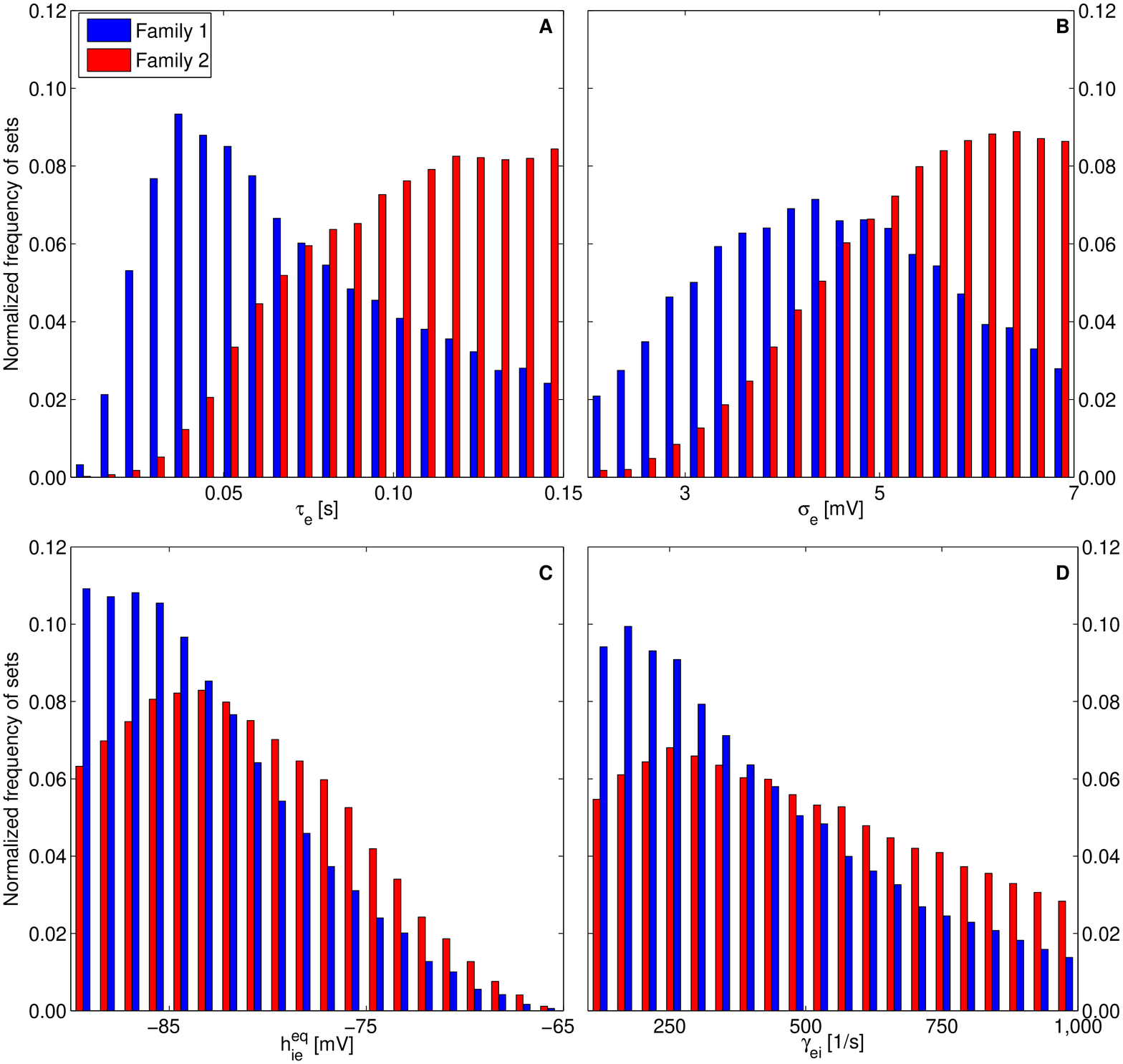}
\end{center}
\caption{{\bf Normalized frequencies of parameter values by families, part 1.} 
Blue bars are for F1, red ones for F2. The total number of bins is $20$, within the intervals 
of physiological validity indicated in Tab.~\ref{tab:allparms}.  
Histograms {\bf A} and {\bf B} represent the distributions 
with the most striking differences.}
\label{fig:distros1}
\end{figure}

The eight largest dissimilarities between the frequency distributions of the
families are found for the parameters $\tau_e~(0.42)$, $\sigma_e~(0.33)$,
$h_{ie}^{eq}~(0.16)$, $\gamma_{ei}~(0.15)$, $\gamma_{ii}~(0.15)$, $\Gamma_{ie}~(0.12)$,
$\Gamma_{ii}~(0.11)$ and $N_{ie}^\beta~(0.11)$, where the number in brackets
is the corresponding $d_{\mathrm{IR}}^j$ value. These distributions are shown 
in Figs.~\ref{fig:distros1} and \ref{fig:distros2}.
In the case of $\tau_e$, see Fig.~\ref{fig:distros1}A, F1 reaches a maximum
around $0.03~\mathrm{s}$, with a slow decline for higher values, while F2 has a steadily increasing
trend, perhaps saturating above $0.1~\mathrm{s}$. This is also the case
for F2 sets in $\sigma_e$, see Fig.~\ref{fig:distros1}B, whereas F1 gives rise to a broad unimodal
distribution with
a maximum around $4.2$ mV. For both $\tau_e$ and $\sigma_e$, the difference in frequency 
for F1 and F2 at the edges of the intervals is striking. This strongly suggests that
$\tau_e$ and $\sigma_e$ play a relevant role in selecting the global dynamics of the model, 
with a clear bias for one type of family over the other.

\begin{figure}[htbf]
\begin{center}
\includegraphics[width=\linewidth]{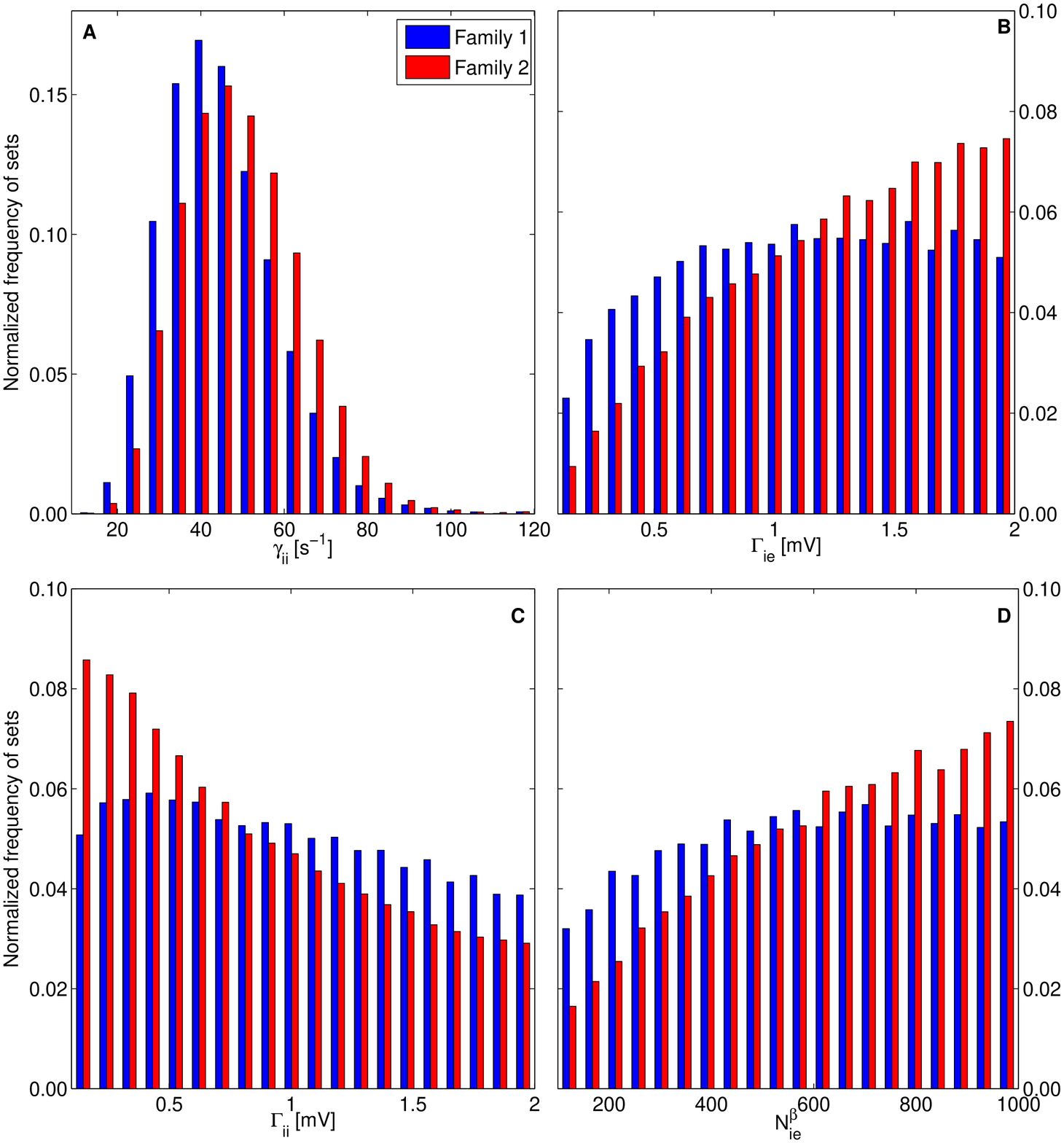}
\end{center}
\caption{{\bf Normalized frequencies of parameters values by families, part 2.} 
See caption of Fig.~\ref{fig:distros1} for details. $\gamma_{ii}$
is here binned only in the interval $0-120/\mathrm{s}$, since larger values are extremely rare.
Also note the different $y\/$-axis scales.}
\label{fig:distros2}
\end{figure}

Although differences in the distributions of other parameters are still significant, they do have  
a minor effect in predicting family membership, since the information radius distance of $\tau_e$ and $\sigma_e$
is by far the largest of all the distributions. For example, 
$h^{eq}_{ie}$ in Fig.~\ref{fig:distros1}C and $\gamma_{ei}$ in Fig.~\ref{fig:distros1}D 
retain separated maxima
for the occurrence of F1 and F2, showing a tendency for larger likelihood of F1 
at low physiological values and a faster decrease in frequency as the parameters increase.
Notice how $\gamma_{ii}$ in Fig.~\ref{fig:distros2}A shows bell-shaped distributions concentrated
in the subinterval $10-100/\mathrm{s}$, in a similar configuration to that of total power
ratios discussed above in Fig.~\ref{fig:allrel}.
In contrast, for PSP amplitudes $\Gamma_{ie}$, $\Gamma_{ii}$, and
connection strengths $N^{\beta}_{ii}$, see Fig.~\ref{fig:distros2}B-D, we observe larger
changes for F2, while F1 is closer to constant throughout the allowed parameter range. 
In these three cases at the edges of the physiological parameter interval one finds clear
bias for one family over the other, but overall the difference is not as strong as for
$\tau_e$ and $\sigma_e$.
The distributions of the remaining parameters do not show marked differences between
F1 and F2, with $d_{\mathrm{IR}}^j$ values between 0.012 and 0.10.

\subsection*{Different oscillatory responses of families to altered inhibition}

Belonging to F1 or
F2 only weakly predicts the presence of a surge or decrease in power when
isoflurane action is simulated, since only the tails in Fig.~\ref{fig:allrel}
have a predominant presence of one family over the other.
Nevertheless, the variabilities of patterns belonging to F1 and F2
embody different responses to changes in inhibition. 
Thus we focus now on the oscillatory activity
associated with F1 and F2, when changes in inhibition are 
modeled solely via alteration of the parameter $R$,
to see how alterations in inhibitory PSP modify the 
oscillatory ``landscape'' families can produce.
To this end we inspected recurrent \textsf{1par} plots at fixed $k$, within a physiological 
range of $0.75 \leq k \leq 1.25$, and characterize the oscillatory activity.

\begin{figure}
\begin{center}
\includegraphics[width=\linewidth]{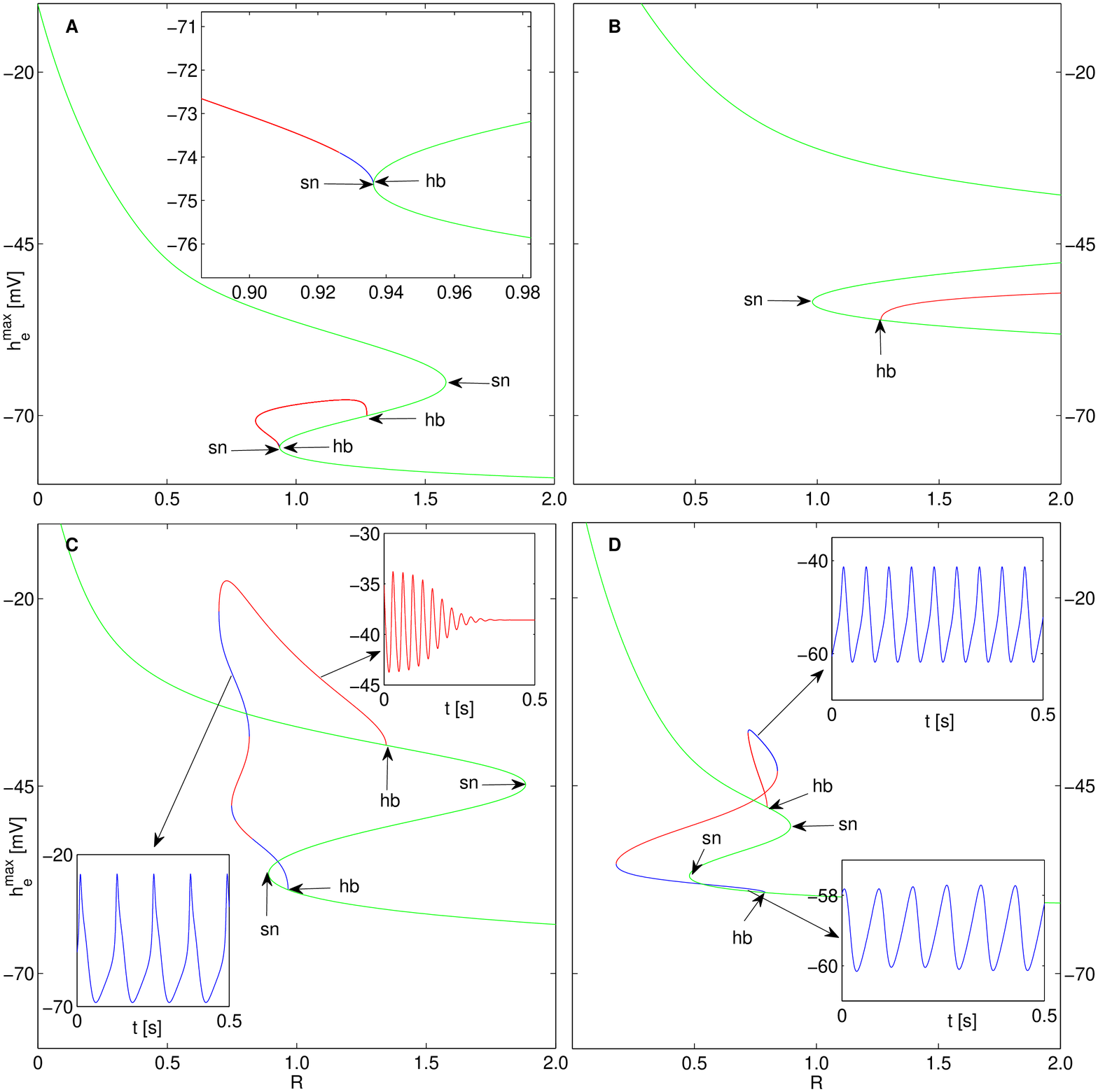}
\end{center}
\caption{{\bf Continuations in $R$ of four different F1 parameter sets
at fixed values of $k$.} 
The interval $0.75 \leq R \leq 2.00$ is considered physiologically 
meaningful, as are the chosen values of $k$
({\bf A} $k = 1.25$, {\bf B} $k = 1.25$, {\bf C} $k = 1.00$ and {\bf D} $k = 1.00$). 
Green lines denote equilibria, blue lines correspond to stable and
red to unstable periodic orbits, i.e., self-sustained oscillations. 
To avoid clutter, only the maximal $h_e$ of orbits are plotted and changes to
the stability of equilibria between the \textsf{sn}s and \textsf{hb}s are not shown. 
Examples of periodic orbits from stable and unstable branches in physiologically 
meaningful intervals are also depicted as insets: in {\bf C} an unstable oscillation converging to equilibrium in red,
and a stable one in alpha frequency range ($8-13~\mathrm{Hz }$) in blue; in {\bf D}
two stable orbits with different amplitudes, the top one in beta ($13-30~\mathrm{Hz }$)
and the bottom one in alpha frequency range, respectively.
Note that the latter two orbits partly overlap in $R$, 
and thus constitute an example of multistability.
}
\label{fig:1par1}
\end{figure}

Examples of this are depicted in Figs.~\ref{fig:1par1} and \ref{fig:1par2},
obtained from selected sets in the initial batch of 405.
We want to stress that the variations within each family are extensive:
the plots proposed are only schematic indications and are not meant to exhaustively describe 
the behaviors of F1 or F2 types. Nonetheless, some general conclusions can be drawn,
based on the most commonly appearing dynamical traits in the continuations in $R$.
Firstly, F1 seems to be less prone to stable oscillatory activity in the 
physiological range than F2. Changes of R, and thus of $\Gamma_{ie}$ and $\Gamma_{ii}$, tend
to trigger unstable periodic orbits with a limited parameter span for F1,
whereas F2 shows stable cycles in a larger parameter range. In particular, plots  
where the size of the
parameter interval of stable orbits is either very limited 
(Fig.~\ref{fig:1par1}A) or zero (Fig.~\ref{fig:1par1}B)
are frequently observed in F1 types, whereas plots with significant oscillatory regimes
(Fig.~\ref{fig:1par1}C and D) are infrequent. 

Secondly, the presence of two or more stable regimes is rarer in sets
belonging to F1, as compared to F2, and furthermore  they typically
extend over more limited  $R$ ranges
within the physiological interval. Hence the occurrence of 
separate periodic orbits at the {\em same} value of $R$ is much rarer
in F1 than in F2. In other words, multistability is distinctive of F2. For instance,  
the plots in Fig.~\ref{fig:1par1}C and D have a very
limited extension of stable orbits for physiological values
of $R$ compared with Figure \ref{fig:1par2}B, C and D.
Since a relation between multistability and memory has been hypothesized, 
this could suggest that F2 types possess
dynamical structure that facilitates the formation of memory.
In general and on average, it appears that F1 
is associated with a more restricted dynamical repertoire than that of F2, and that
F2 responds to the changes in inhibition with more activity than F1.

The distribution of the Hopf bifurcations in different families gives
some insight into this. In F1, the curve of Hopf bifurcations can simply be closed, with no interaction
with the saddle-node bifurcation. This is responsible for the dynamics
in Fig.~\ref{fig:1par1}A and B.
In this case, there exists a time-periodic solution
but it is stable only on a small domain. A slightly more complicated situation
arises if the Hopf bifurcation intersects with the saddle-node bifurcation at
the \textsf{fh} points, as in Fig.~\ref{fig:fams}A. This leads to a more folded branch of periodic solutions,
as shown in Fig.~\ref{fig:1par1}C and D. Again, stable oscillatory behavior is
observed only in a small domain. In F2 instead the \textsf{hb} curve is stretched out and
interacts with both branches of \textsf{sn} curves and generally with both sides of the
$\Lambda$-shaped branch, as in Fig.~\ref{fig:fams}B. This last interaction is, in particular,
closer to the physiological range than the \textsf{fh} points usually observed for F1. 
Hence, the configuration in F2 does have a greater potential for sustained oscillations and the coexistence of
stable periodic solutions, as shown in Fig.~\ref{fig:1par2}B.

\section{Discussion}

We have illustrated how bifurcation diagrams  of Liley's MFM in two 
parameters, $R$ and $k$, account for core effects
of inhibition over the cortex and can be classified into two
topologically distinct families. These parameters modify inhibitory PSPs and
inhibitory self-coupling, respectively, and the two archetypal families 
are identifiable by invariant dynamical features of model cortex. 
Across a large number of 
parameter sets, we have listed and identified all the possible 
diagrams.
The relatively small number of bifurcations and combinations between
branches we have found indicates that Liley's MFM, when it
supports ``realistic alpha activity'' at nominal parameter
values \cite{ref:BL05},  results in a relatively limited dynamical repertoire for
cortical activity in response to parametric variations within physiologically admissible space. 
These dynamical structures
have been classified into families according to their shared global topology
and equivalent local properties of the solutions around
their bifurcations. Family membership weakly correlates with typical changes in 
total spectral power when the cortex is acted upon by 
anesthetic agents, but strongly with unusually large changes. At the same time,
the likelihood for family membership is sensitive to changes in specific neurophysiological attributes,
in particular to the standard deviation of the excitatory neural population firing threshold $\sigma_e$ (i.e.
the steepness of the neuronal firing rate function $S_{k}$) and
the excitatory membrane decay time constant $\tau_e$. As far as we are aware, 
these connections between dynamical properties of a MFM,
physiological attributes of cortex and changes in
EEG power spectra under anesthesia have never been found before.

\begin{figure}
\begin{center}
\includegraphics[width=\linewidth]{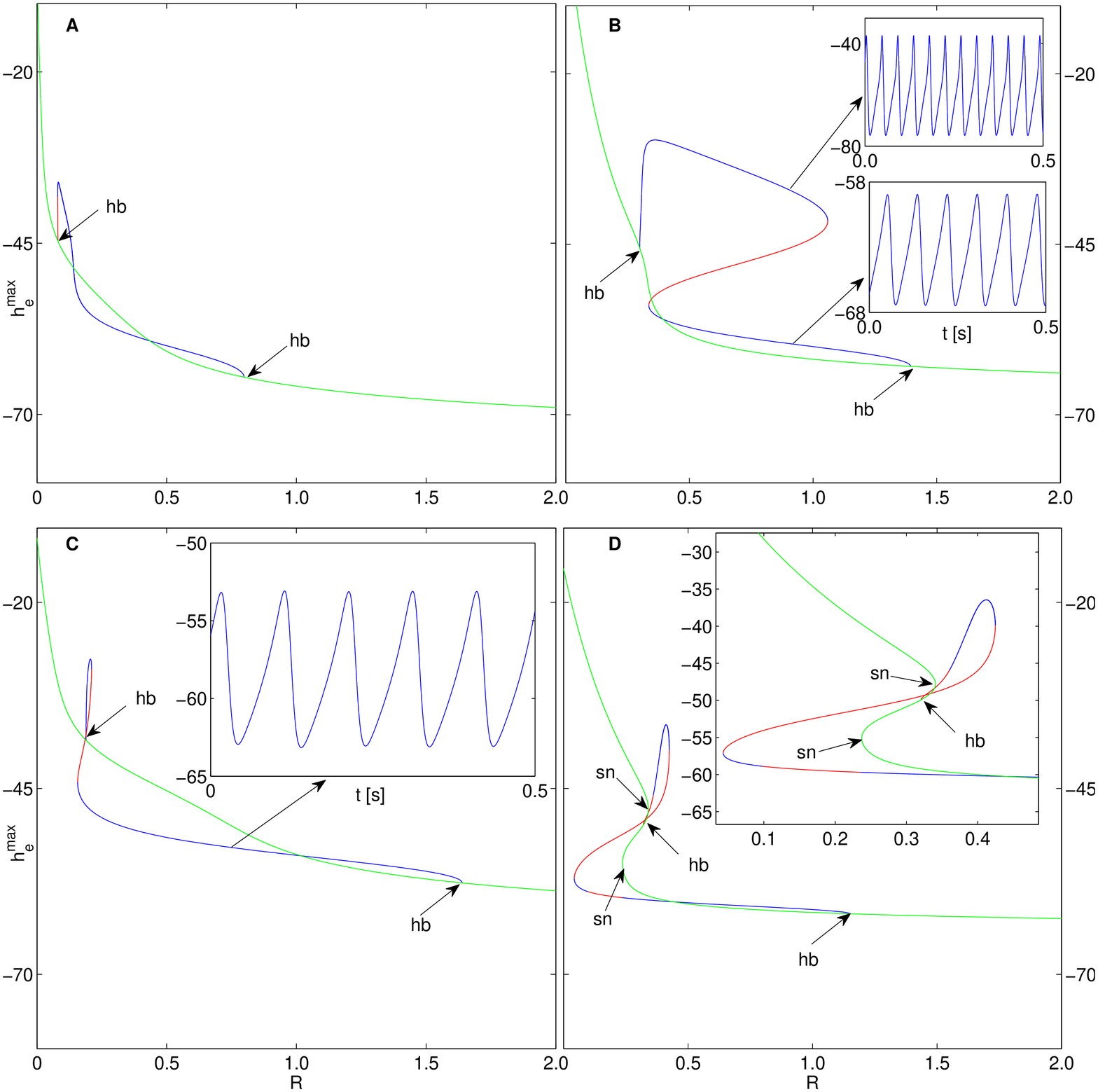}
\end{center}
\caption{{\bf Continuations in $R$ of four different F2 parameter sets
at fixed values of $k$.} 
As in Fig.~\ref{fig:1par1}, but for F2 and with the following values of $k$: 
{\bf A}  $k = 1.00$, {\bf B} $k = 1.25$, {\bf C} $k = 1.25$ and {\bf D} $k = 1.25$.
Note how the extension
of stable orbits and the existence of multistable regimes (i.e., multiple,
different orbits occurring at the same $R$) is more pronounced compared to F1
in Fig.~\ref{fig:1par1}. In {\bf B} and {\bf C} stable periodic orbits at physiologically meaningful values are also shown
as insets. In {\bf B}, where a multistable regime is present, frequencies are in the beta and alpha bands for the higher and lower inset, respectively, whereas in
{\bf C} a typical alpha oscillation is depicted. The inset in {\bf D} is a magnification showing more clearly the presence of multistability in the non physiological $R$ range.
}
\label{fig:1par2}
\end{figure}

It is important to appreciate the generality and relative simplicity of
the methods we have used here. We have chosen two meaningful parameters for the bifurcation analysis, $R$ and $k$,
and have then classified bifurcation diagrams into
two families, F1 and F2, according to the absence or presence of
separate \textsf{sn} lines. Family membership turned out to be separated
by a swallow tail bifurcation. More sophisticated criteria for labeling resulting
\textsf{2par} plots could be envisaged. For example, families might have been
categorized by the number, location and types of bifurcation they presented, or by the values
of the coefficients in their normal forms. However, in our case the 
simple approach proved sufficient. The main methodological result we wish to underline 
is that by grouping the global dynamics 
for a representative sample of 73,454 parameter sets
of an EEG model according to
agreed criteria, strong, distinctive and robust
correlations with EEG
power spectra on one hand and some of the physiological parameters on the other hand
have been derived.
This is novel {\it per se}, because it shows that systematic
partitioning of complex parameter spaces according to global dynamical patterns can lead
to significant results which are not otherwise accessible. This method seems to be one the few 
viable strategy for unveiling interesting relations in high dimensional nonlinear spaces and
exploring typical responses in complex models.

In this paper, we have focused on the action of anesthetics and ratios of total power, but 
many other developments are possible. 
For example, we could select those sets which give rise to epileptic
activity or have a high power in the EEG gamma band, categorize their
dynamical patterns and look for parameter correlations.
Our approach is hence promising for
other theoretical studies of the EEG, or neural activity in general.

We want to stress that the proposed method is not dependent on Liley's MFM: in general,
any other set of coupled, nonlinear ODEs with a complicated but meaningful parameter space 
can be treated likewise. The difficulties in extracting
relevant information from equations with many parameters 
and involved feedback and feedforward mechanisms should not be underestimated,
and our dynamical approach can offer new guidance. Here
a PCA on the 73,453 analyzed parameter sets turns out to be
completely useless. This is also true if we employ a PCA after dividing the sets 
into F1 and F2: the first ten principal components then account for about 55\% of
the total variance, up from 42\%. Again, no valuable information on parameter relations can
be gleaned from these
types of statistical methods. If we compare Fig.~\ref{fig:pareto}
with the clear separation between distributions in power ratio in Fig.~\ref{fig:allrel}
or in Fig.~\ref{fig:distros1} for parameters $\tau_e$ and $\sigma_e$, 
we can appreciate the strength of our method to tease out hidden correlations.
 
What is the biological meaning of family differences with respect to
parameters such as $\tau_e$ and $\sigma_e$? 
First, a small excitatory rate constant $\tau_e$ enhances the speed with which
$h_e(t)$ adapts to changes in PSP inputs, as is clear from Eq.~\ref{eq:liley1}.
For $\tau_e=0$, the mean excitatory soma membrane potential simply becomes
a weighted sum of the PSPs. Thus, low $\tau_e$ enslave $h_e(t)$ more directly
to synaptic input and results in a depression of self-sustained oscillations, which
is characteristic for F1. $\sigma_e$ parametrizes the local excitatory firing rate
response, see Eq.~\ref{eq:liley3}, and hence does not directly act upon $h_{e}(t)$
dynamics, as $\tau_e$ does.  However, a smaller $\sigma_e$ means that the
excitatory firing rate changes more rapidly for variations of $h_e(t)$, although over
a more limited range. In the limit
of  $\sigma_e=0$, one obtains a step function which instantaneously switches the population
from zero to maximal firing at  the threshold value $\mu_e$. As mentioned above,
in a steady state situation excitatory self-feedback is essentially proportional
to $(N^\alpha_{ee}+N^\beta_{ee})\cdot S_e\left[h_e(t)\right]$, and a low $\sigma_e$
increases this self-feedback. Thus, as for $\tau_e$,  the reaction of $h_e(t)$ to changing synaptic
input will be more rapid for low $\sigma_e$, limiting the presence of self-sustained oscillations
and hence favoring membership in F1. One can then view the somewhat lesser differentiation of
$\sigma_e$ with regards to the families due to its boosting excitatory reactivity
more indirectly. Overall, enhancements in the reactivity of excitatory populations result in a diminished
dynamical repertoire.

These considerations on the primary role of $\sigma_e$ and $\tau_e$ in influencing
family membership are not weakened by the presence of complex feedbacks between neuronal populations
in Liley's model, c.f. Fig.~\ref{fig:scheme}. Sets of F1 types appearing
at low $\tau_e$ or $\sigma_e$ show homogeneous distributions for all
the other parameters, with no formation of clusters or any other form of bias. This suggests that
the model cortex generally can be pushed towards F1 by appropriately reducing the excitatory 
membrane decay constant or the excitatory standard deviation of the firing threshold, and towards F2 by
increasing them. Thus, we can also identify excitatory reactivity as an endogenous counterpart to the
exogenous control.
In fact, it has been shown in detail here also how exogenous activity, namely the excitatory thalamic 
input to inhibitory (excitatory) cortical populations described by $p_{ei}$ ($p_{ee}$), 
can cause changes across and within families.
There exists a kind of mirror symmetry between $p_{ee}$ and $p_{ei}$: at 
the increase of $p_{ei}$, or decrease of $p_{ee}$, 
we have observed a transition from F2 to F1, see Fig.~\ref{fig:F1toF2}, corresponding to 
a transition from complex to simple dynamics. This indicates
that the exogenous driving of inhibition (or a reduction in the driving of excitation)
moves the cortex to dynamics which, on average, appear simpler.
  
Thalamic input appears in Liley's MFM as an active source of control
on the complexity of brain dynamics, having the role of modulator in
cortical interactions.  This is quite different to the classical view
of thalamus as being the gateway through
which peripherally derived sensory information reaches cortex\cite{ref:W38}. For
example, retinal axons project to the lateral geniculate nucleus of the
thalamus, which then projects to visual cortex.  In this way the
thalamus is seen as the conduit through which all information about lower
levels of the nervous system reaches the cortex for ``processing''.
Viewed from the perspective of our model, this would imply that
thalamic input simply drives the cortex.  However, based on our 
metabifurcation analysis we predict that such input also configures (i.e.
``modulates'') local and global dynamical responses from the cortex.

It is now known that all areas of cortex, not just the primary sensory
cortices (somatosensory, visual, auditory) receive thalamic input, and
that most areas of cortex, in some reciprocal manner, innervate
thalamus\cite{ref:SG06}.  As a consequence, contemporary approaches to understanding
the role of thalamus are now aimed towards better articulating the
functional inter-relationships that exist between cortex and thalamus.
One popular approach has been the
modeling of integrated dynamics of thalamo-cortical activity.  On this
basis, it has been speculated that the human alpha rhythm emerges as a
consequence of reverberant activity between cortex and thalamus.  While
we have not here considered such thalamo-cortical feedback, doing so in
the context of our results leads to a number of interesting
speculations regarding the auto-regulation of cortical dynamics:
cortical feedback through the thalamus could initiate a sequence of
transitions between bifurcation pattern families, providing a way of
reconfiguring cortical dynamics ``on the fly'' and thus on a time scale
quite different to that associated with activity dependent synaptic
plasticity. Unfortunately, not enough is known experimentally about
thalamocortical feedback at present to turn such speculations into
strong model predictions.

On the other hand, thalamus also appears to stimulate internal variations
{\em within} families, since a lower value of $p_{ei}$ (higher value of $p_{ee}$)
corresponds to more complex bifurcations as depicted in Fig.~\ref{fig:AnnBTs1}. 
As said, this type of \textsf{2par} plot with two \textsf{bt} points 
has been shown to produce chaotic
activity, meaning that ``cortical state'' F2
is capable of producing highly nonlinear events, and an
exogenous enhancement of cortical inhibition tends to reduce this ability.
This and the fact that F2 types have shown to occur with the largest likelihood in parameter space (i.e.,
their total number is one and a half time more than F1), 
stimulates some speculations. F2 may be considered as the cortical default state
of Liley's MFM, associated with the richest dynamical responses to inhibition
and acting as a sort of cognitive readiness state. Both endogenous (cortical) enhancements of
excitatory reactivity and exogenous, subcortical increases (decreases) of inhibitory (excitatory) activity reduce
the cortical state to much simpler behavior. 

Finally, we have implemented a partition of the parameter space
according to the bifurcation plots. As mentioned, 
types of normal forms appear to be constant within each family 
and continuation lines of the diagrams are similar, 
whilst the extensions of these patterns in the $Rk$-space vary wildly. 
This suggest the hypothesis that a single, high dimensional master bifurcation 
diagram exists and governs all sets, and that
changes in parameters only change its orientation, position and scale in the (hyper)volume. By
continuing in $R$ and $k$, we intersect the master diagram with a plane.
The position of this plane determines whether we see F1 or F2
diagrams, and the two can be continuously deformed into each other.

We would like to conclude with some future directions. First, an inventory
of all the possible types of oscillatory dynamics associated with F1 and F2
is surely attractive. How do many of them support complex events such as chaos and how is this
reflected in parameter space? Do chaos and complex dynamics 
occur for limited, selected intervals in the parameters?
The great hurdle is represented by the computational
demand: it is still unfeasible to process all sets one by one and see what types of orbits are present.
Second, a systematic study of the unfolding of the bifurcations in \textsf{2par} plots that
characterize families should be attempted. This has the potential to increase our knowledge about 
the type of dynamics these parameter sets can produce.
Finally, there is an interesting theoretical question concerning the optimization of the
method we presented in this paper. What is the best choice of continuation parameters 
for achieving the highest degree of separation among bifurcation plots? 
We wonder if it would be possible to choose parameters or combinations of parameters such that their
degree of clustering with respect to families is maximal. Some canonical 
discriminant analysis on the sets has been preliminarily performed, but results are not
particularly illuminating. Achieving optimal partitioning has the potential to uncover biologically relevant dependencies among
parameters that are still unknown, and to shed more light on 
clinically or pharmacologically relevant questions.

\section{Acknowledgements}

FF was supported by Australian Research Council Discovery Grant
DP 0879137 and the Swinburne University of Technology Researcher
Development Scheme 2010. FF would like to acknowledge the hospitality of
the Department of Mathematics at Concordia University, Canada, where
the idea for this study was conceived.  LvV was supported by NSERC
Discovery Grant 355849-2008.  IB would like to thank Prof.\ Rolf
K{\"{o}}tter for the time granted to complete this project.

\bibliographystyle{elsart-num-sort}
\bibliography{biblio}

\end{document}